# Successive Estimations of Bilateral Bounds and Trapping/Stability Regions of Solution to Some Nonlinear Nonautonomous Systems

Mark A. Pinsky

**Abstract.** Estimation of the degree of stability and the bounds of solutions to nonautonomous nonlinear systems present major concerns in numerous applied problems. Yet, current techniques are frequently yield overconservative conditions which are unable to effectively gage these characteristics in time-varying nonlinear systems. This paper develops a novel methodology providing successive approximations to solutions that are stemmed from the trapping/stability regions of these systems and estimate the errors of such approximations. In turn, this leads to successive approximations of both the bilateral bounds of solutions and the boundaries of trapping/stability regions of the underlying systems. Along these lines we formulate enhanced stability/boundedness criteria and contrast our inferences with inclusive simulations which reveal dependence of the trapping/stability regions upon the structure of time – dependent components and initial time- moment.

1. **Introduction**

Stability and boundedness of solutions to nonlinear systems with variable coefficients are of major concerns in various applications. For instance, the corresponding stability conditions aid the design of robust controllers and observers [1-3]. Sufficient conditions of stability of fixed solutions to such systems were initially established by Lyapunov [4] and subsequently advanced in various publications, see [5-7] and review paper [8] for additional references. In control literature, the method of Lyapunov functions has been used extensively for stability analysis of some nonlinear and nonautonomous systems [9-15]. Yet, adequate Lyapunov functions are rare in this field. Convenient conditions embracing boundedness/stability of solutions to nonhomogeneous nonlinear systems with variable coefficients are practically unknown as well.

Nonetheless, local stability conditions are well known for nonlinear systems with time – invariant and Hurwitz linearization and smooth or Lipschitz continuous nonlinear perturbations [1, 2, 5]. Such conditions can be naturally extended for nonlinear systems with periodic linearization. Yet, local dynamics of systems with nonperiodic linearization turn out to be more complex and sensitive to small perturbations. Consequently, assessment of local stability of such systems requires specialized techniques resting on application of Lyapunov and generalized exponents [4, 5].

Let us briefly review some explicit stability conditions for nonautonomous nonlinear systems which can be written as follows,

$$\dot{x} = A(t)x + f(t,x) + F(t), \quad \forall t \geq t_0, \quad x \in H \subset \mathbb{R}^n, \quad f(t,0) = 0$$
$$x(t_0, t_0, x_0) = x_0 \tag{1.1}$$

where, $f := [t_0, \infty) \times \mathbb{R}^n \to \mathbb{R}^n$ and $F := [t_0, \infty) \to \mathbb{R}^n$ are continuous functions, matrix, $A(t) \in \mathbb{R}^{n \times n}$ is continuous, $t \geq t_0$, $t_0 \in T := [\varsigma, \infty)$, $\varsigma \in \mathbb{R}$, $x_0 \in \vartheta \subset \mathbb{R}^n$, $\vartheta$ is a bounded neighborhood about $x \equiv 0$, $F(t) = F_0 \eta(t)$, $\sup_{t \geq t_0} \|\eta(t)\| = 1$, $F_0 \in \mathbb{R}_{\geq 0}$, $\mathbb{R}_{\geq 0}$ is a set of all nonnegative real numbers, $\|\cdot\|$ stands for induced 2-norm of a matrix or 2-norm of a vector and $x(t, t_0, x_0) := [t_0, \infty) \times T \times \mathbb{R}^n \to \mathbb{R}^n$ is the solution to (1.1). To shorten the notation, we will write below that $x(t, t_0, x_0) \equiv x(t, x_0)$ and assume that (1.1) possesses a unique solution $\forall t \geq t_0$.

We will also reference the homogeneous counterpart of (1.1),

$$\dot{x} = A(t)x + f(t,x)$$
$$x(t_0, x_0) = x_0 \tag{1.2}$$

and its linearization,

$$\dot{x} = A(t)x$$
$$x(t_0, x_0) = x_0 \tag{1.3}$$

Obviously, the solution to (1.3) can be written as $x(t, x_0) = W(t, t_0) x_0$, where $W(t, t_0) = w(t, t_0) w^{-1}(t_0, t_0)$ and $w(t, t_0)$ are transition and fundamental solution matrices for (1.3), respectively.

Mark A. Pinsky. Department of Mathematics and Statistics, University of Nevada.Reno, Reno NV 89557, USA.
e-mial:pinsky@unr.edu. **To appear in Journal of Nonlinear Dynamics.**

Contemporary stability conditions of the trivial solution to (1.2) are frequently derived under Lipschitz continuity condition,

$$\|f(t,x)\| \le l_1(t,t_0)\|x\|, \quad \forall x \in \Omega_1 \subset \mathbb{R}^n, \forall t \ge t_0 \tag{1.4}$$

(where $\Omega_1$ is a bounded neighborhood of $x \equiv 0$ and $l_1(t,t_0) \le \hat{l}_1 \in \mathbb{R}_{>0}, \forall t \ge t_0$ is a continuous function) and some bounds on the norm of transition matrix [5- 8],

$$\|W(t,t_0)\| \le N e^{-v(t-t_0)}, \quad \forall t \ge t_0, N, v > 0 \tag{1.5}$$

or

$$\|W(t,t_0)\| \le N \exp \int_{t_0}^{t} \Upsilon(s) ds, \forall t > t_0 > 0 \tag{1.6}$$

Then, asymptotic stability of the trivial solution to (1.2) is embraced if either (1.4), (1.5) and

$$N\hat{l}_1 - v < 0 \tag{1.7}$$

hold [5, 6] or (1.4), (1.6) and

$$\limsup_{t \to \infty} (1/(t-t_0)) \int_{t_0}^{t} \Upsilon(s) ds + N\hat{l}_1 < 0 \tag{1.8}$$

hold [7, 8]. Clearly, (1.6) reduces to (1.5) if $\Upsilon(t) = const$, otherwise, to our knowledge, this function was not explicitly defined yet for systems of practical concerns. Utility of Lipschitz constant, i.e. $\hat{l}_1$, in (1.7) and (1.8) usually contributes to overconservative nature of these conditions [16, 17]. Additionally, efficient bounding of $\|w(t,t_0)\|$ can present a challenging task even if $w(t,t_0) = \exp A(t-t_0)$ with $A = const$ [18]. More efficient stability conditions of trivial solution to (1. 2) were developed in [19], but verification of these conditions can be challenging by its own.

   Virtually all known stability conditions for nonlinear time-varying systems turn out to be overconservative and local in nature - they unable to adequately estimate the region of stability for these systems. However, the problem of estimating the stability regions for autonomous nonlinear systems has attracted a flood of publications in the past few decades, see recent reviews in [20,21]. Most of the efforts in this area are based on optimizing the Lyapunov functions and combining this methodology with the level set method [20 – 28]. Another popular tactic is trajectory reversing [29 – 33] which is also combined with level set method [34, 35]. There are also other techniques exploring geometrical reasoning [36-38] and some other strategies [39]. Yet, neither of these techniques were targeted nonautonomous nonlinear systems at the present time.

   So-known geometric methods were developed to stability analysis of nonlinear and nonlinear controlled systems, see, e.g., [2, 3, 49] and subsequent references therein.

   Though stability assesses the fate of solutions on long time – intervals, estimation of their evolutions on initial/intermediate time – intervals is frequently equally imperative [40 – 42]. Nonetheless, this problem has not received sufficient attention in the contemporary literature.

   In [43] under some conditions we develop an approach to estimation of the upper bound of norms of solutions to (1.1)/(1.2) which leads to more general boundedness/ stability criteria and estimation of the trapping/stability regions for these systems. Still, the developed estimates turn out to be rather conservative, especially for nonhomogeneous systems and ones with large nonlinear components.

   This paper develops a novel methodology which successively approximates the solutions stemming from trapping/stability regions of systems (1.1)/(1.2) and estimate the norms of errors of these approximations. This leads to recursive estimation of both the bilateral bounds of solutions' norms and the trapping/stability regions. Along these lines we formulate novel bondedness/stability criteria and reinforce our theoretical results in representative simulations. Our inferences advance the methodology developed in [43] and uplift some of its major limitations.

   This paper represents a solution of (1.1) as,

$$x(t, x_0) = Y_m(t, x_0) + z(t, x_0), \ 1 \leq m \subset \mathbb{N}, \ t \geq t_0 \tag{1.9}$$

where $Y_m := [t_0, \infty) \times \mathbb{R}^n \to \mathbb{R}^n$, $x(0, x_0) = Y_m(0, x_0) = x_0$, $z(0, x_0) = 0$ and $\mathbb{N}$ is a set of positive integers. Thus, $Y_m(t, x_0)$ can be viewed as an approximate solution to (1.1) assuming its initial data and $z(t)$ is considered as the error-function of this approximations that depends implicitly upon $x_0$. Clearly, (1.9) yields the following bilateral bounds,

$$\|Y_m(t, x_0)\| - \|z(t, x_0)\| \leq \|x(t, x_0)\| \leq \|Y_m(t, x_0)\| + \|z(t, x_0)\|, \ t \geq t_0 \tag{1.10}$$

where the left bound should be adjusted to zero if it takes negative values.

Next, we represent $Y_m(t, x_0)$ in the following form,

$$Y_m(t, x_0) = \sum_{i=1}^{m} y_i(t, x_0) \tag{1.11}$$

where $y_i := [t_0, \infty) \times \mathbb{R}^n \to \mathbb{R}^n$ are solutions to some linear equations which are derived subsequently.

Afterward, we show under some conditions that either $\lim_{t \to \infty} \|Y_m(t, x_0)\| = 0$ if $F_0 = 0$ and $x_0$ belongs to the stability region of the trivial solution of (1.2) or $\lim_{t \to \infty} \|Y_m(t, x_0)\| \leq O(F_0)$ if $F_0 > 0$ and $x_0$ belongs to the trapping region of (1.1). Thus, in both cases the behavior of $\|x(t, x_0)\|$ on long time intervals is determined by evolution of $\|z(t, x_0)\|$.

This paper is organized as follows. Next section presents the synopsis of some essential technical results. Sections 3 and 4 derive equations for $y_i(t, x_0)$ and $z(t, x_0)$ for a relatively simple system and, subsequently, in a general form. Section 5 examines behavior of the approximate solutions on long time – intervals, section 6 derives a linear scalar equation for the upper bond of $\|z(t, x_0)\|$ and develops its applications, section 7 derives a nonlinear counterpart of this equation and apply it for estimating of the bilateral bounds of solutions and the trapping/stability regions. Section 8 review the results of our simulations and section 9 concludes this study and outlines the topics of subsequent works. The appendix derives some technical results pertained to section 7.

## 2. Preliminaries

### 2.1. Lyapunov and General Exponents

Firstly, let us briefly contrast the definitions of broadly used the Lyapunov exponents, see, e.g. [8] and the general exponents, which were introduced in [5]. The latter frequently provide more conservative inferences but possess "more natural properties". In fact, the spectrum of Lyapunov exponents for (1.3) sensitively depends upon perturbations of the underlying system [44-46] whereas the maximal/minimal general exponents are robust under such perturbations [5] which enables our inferences in section 5. Different terminology and notations are currently used for these measures and we adopt ones from [8] and [5], respectively.

For equation (1.2) the Lyapunov exponents, $\chi_i(x(t, x_0)) = \lim_{t_0 < t \to \infty} \sup (t - t_0)^{-1} \ln \|x(t, x_0)\|$, $i \leq n$, where $x(t, x_0)$ is a solution to (1.2). The maximal Lyapunov exponent, $\chi = \max_i \{\chi_i(x(t, x_0))\}$ measures the maximal rate of exponential growth/decay of the corresponding solutions as $t \to \infty$ and bears a pivoting role in stability theory. For a linear system (1.3) the Lyapunov exponents are defined as, $\mu_i = \lim_{t_0 < t \to \infty} \sup (t - t_0)^{-1} \ln \sigma_i(\|w(t)\|)$, $i \leq n$, where $\sigma_i$ are the singular values of the fundamental matrix of (1.3); and for linear systems, $\max_i \mu_i = \chi$ [8]. This let to bound the transition matrix of (1.3) as follows,

$$\|W(t,t_0)\| \leq D_\varepsilon \exp((\chi+\varepsilon)(t-t_0)), \ D_\varepsilon > 0, \ \forall t \geq t_0 \tag{2.1}$$

where $\varepsilon$ is a small positive number [5].

In turn, the upper general exponent, $\kappa$, is defined in [5] as follows. Let $x(t,x_0)$ with $\|x_0\| \leq r$ be a solution of (1.2) stemming from a centered at zero ball with radius $r$. Assume that every of such solutions obeys the inequality

$$\|x(t,x_0)\| \leq N \exp(-\nu(t-t_0))\|x(t_0)\|, \ N,\nu > 0, \ \forall t \geq t_0 \tag{2.2}$$

Then, infimum of values of $\kappa = -\nu$, for which (2.2) holds, is called the upper general exponent for (1.2) at zero. Clearly, $-\kappa = \sup \nu$ for which (2.1) can be defined. It is stated in [5] that $\chi \leq \kappa$ as well as that the trivial solution to (1.2) is exponentially stable if $\kappa < 0$.

For linear equation (1.3) condition $\|x_0\| \leq r$ is dropped from the definition of $\kappa \equiv \kappa_A$ and for this equation (2.2) implies that,

$$\|W(t,t_0)\| \leq N \exp(-\nu(t-t_0)), \ N,\nu > 0, \ \forall t \geq t_0 \tag{2.3}$$

Note also that (2.3) is the necessary and sufficient condition for asymptotic stability of (1.3) [2].

In section 5 we use some of the theorems from sections 3.4 and 3.5 in [5] which are abbreviated below.

**Statement 1**. Consider a linear equation,

$$\dot{x} = (A(t) + B(t))x, \ x \in \mathbb{R}^n \tag{2.4}$$

where both, $A(t)$ and $B(t) \in \mathbb{R}^{n \times n}$, $t \geq t_0$ are continuous matrices. Assume that,

$$\lim_{t,s \to \infty} \sup s^{-1} \int_t^{t+s} \|B(\tau)\| d\tau < \delta \tag{2.5}$$

Then, $\forall \varepsilon > 0$ there is such $\delta(\varepsilon)$ that, $\kappa_{A+B} \leq \kappa_A + \varepsilon$, where $\kappa_A$ and $\kappa_{A+B}$ are maximal general exponents for equations (1.3) and (2.4), respectively. Note that (2.5) is contented if, e.g., $\|B(t)\| \leq \delta, \ \forall t \geq t_* \geq t_0$.

Furthermore, it follows from the above statement that if $\kappa_A < 0$ and $\delta$ is sufficiently small, then, $\kappa_{(A+B)} < 0$ as well.

**Statement 2.** Assume that $\delta = 0$. Then, $\kappa_{A+B} = \kappa_A$ and
$\|W_{A+B}(t,t_0)\| \leq N_* \exp(\kappa_A(t-t_0)), \ N_* > 0, \ \forall t \geq t_0$, where $W_{A+B}(t,t_0)$ is a transition matrix of (2.4). The last conditions are intact if, e.g. $\lim_{t \to \infty} \|B(t)\| = 0$.

### 2.2. Comparison Principle

This paper references the comparison principle [2] which is outlined as follows for convenience. Let us consider the IVP for a scalar differential equation,

$$\begin{aligned} \dot{x} &= g(t,x), \ \forall t \geq t_0, \ x(t) \in \mathbb{R} \\ x(t_0) &= x_0 \end{aligned} \tag{2.6}$$

where function $g(t,x)$ is continuous in $t$ and locally Lipschitz in $x$ for $\forall x \in \wp \subset \mathbb{R}$. Suppose that the solution to (2.6), $x(t,x_0) \in X_+, \ \forall t \geq t_0$. Next, consider a differential inequality,

$$\begin{aligned} D^+ y &\leq g(t,y), \ \forall t \geq t_0 \\ y(t_0) &= y_0 \leq x_0 \end{aligned} \tag{2.7}$$

where $D^+ y$ denotes the upper right-hand derivative of solution to (2.6) [2] and $y(t, y_0) \in \wp$, $\forall t \geq t_0$. Then, $x(t, x_0) \geq y(t, y_0)$, $\forall t \geq t_0$.

**2.3. Estimation of Solutions' Norms**

Next, we recap our approach for estimating the norm of solutions to (1.1) which is used in sections 6 and 7, see [43] for additional details. This paper derives a differential inequality for the norms of solutions to (1) in the following form,

$$D^+ \|x\| \leq p(t,t_0)\|x\| + c(t,t_0)\|f(t,x) + F(t)\|$$
$$\|x(t_0, x_0)\| = \|w^{-1}(t_0, t_0) x_0\| = X_0$$
(2.8)

where

$$p(t,t_0) = d\left(\ln \|w(t,t_0)\|\right)/dt$$
(2.9)

and

$$c(t,t_0) = \|w(t,t_0)\| \|w^{-1}(t,t_0)\| = \sigma_{\max}(w)/\sigma_{\min}(w)$$
(2.10)

is the running condition number of $w(t,t_0)$, $\sigma_{\max}(t)$ and $\sigma_{\min}(t)$ are running maximal and minimal singular values of $w(t)$.

Note also that [43] assumes that $\|w(t_0,t_0)\| = 1$ and, hence, (2.9) can be written as follows,

$$\|w(t,t_0)\| = \exp\left(\int_{t_0}^{t} p(s) ds\right)$$
(2.11)

Note that below we adopt a shorter notation: $w(t,t_0) = w(t)$, $p(t,t_0) = p(t)$, $c(t,t_0) = c(t)$ and assume that $c(t) \leq \hat{c} = const$.

Due to comparison principle, solutions to (2.8) are bounded from above by solutions of the following scalar equation,

$$\dot{X} = p(t) X + c(t) \|f(t, x(t, x_0)) + F(t)\|$$
$$X(t_0, X_0) = \|w^{-1}(t_0) x_0\| = X_0$$
(2.12)

Consequently, $X(t, X_0) \geq x(t, x_0)$, $\forall t \geq t_0$.

Furthermore, application of (1.4) and standard norm's inequalities linearizes (2.12) for $\forall x(t, x_0) \in \Omega_1$, $\forall t \geq t_0$ and let us to write this equation as follows,

$$\dot{X}_l = \lambda(t) X_l + c(t)\|F(t)\|$$
$$X_l(t_0, X_0) = \|w^{-1}(t_0) x_0\| = X_0$$
(2.13)

where $\lambda = p + c l_1$. Application of (2.13) enabled us to bound the solutions to (1.1) and derived some boundedness and stability criteria for (1.2). To refine these estimates, we formulated in [43] a nonlinear extension of (1.4), i.e.,

$$\|f(t,x)\| \leq L(t, \|x\|), \forall x \in \Omega_2 \subset \mathbb{R}^n, \forall t \geq t_0$$
(2.14)

where $\Omega_2$ is a bounded neighborhood of $x \equiv 0$ and function $L := [t_0, \infty) \times \mathbb{R}_{\geq 0} \to \mathbb{R}_{\geq 0}$ is continuous in $t$ and $\|x\|$, $L(t, 0) = 0$. Furthermore, [43] defined $L(t, \|x\|)$ explicitly if $f$ is either a polynomial in $x$ or can be approximated by such polynomial with, e.g. Lagrange error term. Let us show for convenience how to define $L(t, \|x\|)$ in a simple but typical case. Assume that $f = \begin{bmatrix} a_1(t) x_1^2 x_2 & a_2(t) x_1 x_2 \end{bmatrix}^T$. Then, $\|f\|_2 \leq \|f\|_1 \leq |a_1| |x_1^2| |x_2| + |a_2| |x_1| |x_2| \leq |a_1| \|x\|^3 + |a_2| \|x\|^2$.

Thus, (2.14) turns into a global inequality for polynomial vector fields or ones which are approximated by polynomials with globally bounded error terms.

Application of (2.14) turns (2.8) into the following inequality,

$$D^+ \|x\| \leq p(t) \|x\| + c(t) \left( L(t, \|x_0\|) + \|F(t)\| \right)$$
$$\|x(t_0, x_0)\| = \|w^{-1}(t_0) x_0\| = X_0 \tag{2.15}$$

which, due to comparison principle, yields the equation,

$$\dot{X}_* = p(t) X_* + c(t) \left( L(t, X_*) + \|F(t)\| \right)$$
$$X_*(t_0, X_0) = \|w^{-1}(t_0) x_0\| = X_0 \tag{2.16}$$

where $X_*(t, x_0) \geq \|x(t, x_0)\|$, $\forall t \geq t_0$. This last equation frequently enables more accurate estimation of the bounds of solutions to (1.1), derivation of enhanced stability criteria and estimation of the trapping/stability regions for equations (1.1) and (1.2), respectively.

Additionally, [43] derives that $\|x(t, x_0)\| \leq c(t_0) \exp\left( \int_{t_0}^{t} p(s) ds \right) \|x_0\|$, where $x(t, x_0)$ is a solution to (1.3).

Thus, the last formula implies that,

$$\kappa_A = -\nu \leq \sup_{t > t_0} \left\{ (t - t_0)^{-1} \int_{t_0}^{t} p(s) ds \right\} \tag{2.17}$$

Next, assume that $\lambda(t, t_0) < -\hat{\lambda}(t_0)$, $\forall t \geq t_0$, $0 < \hat{\lambda}(t_0) \in \mathbb{R}$. Then, $(t - t_0)^{-1} \int_{t_0}^{t} p(s) ds < -\hat{\lambda}$ since $k(t) l_1(t) > 0$ and, in turn, $\nu \geq \hat{\lambda}$. This last inequality is used in section 6 of this paper.

Finally, we extend some typical definitions of trapping/stability regions used for autonomous systems on nonautonomous systems (1.1) and (1.2), respectively. Let us firstly acknowledge a set of initial vectors, $x_0 \in \Im(t_0) \subset \mathbb{R}^n$, $\forall t_0 \in T$, where it is assumed that $0 \in \Im(t_0)$. Then, we present below the following definitions.

Definition 1. A compact set of initial vectors, which includes zero-vector, is called a trapping region of equation (1), if condition, $x_0 \in \Im(t_0)$ implies that $x(t, t_0, x_0) \in \Im(t_0)$, $\forall t \geq t_* \geq t_0$.

Definition 2. An open set of initial vectors, which includes zero-vector, is called a stability region of the trivial solution to (2) if condition, $x_0 \in \Im(t_0)$ implies that $\lim_{t_0 \leq t \to \infty} x(t, t_0, x_0) = 0$.

We illustrate in section 8 that the trapping and stability regions for nonautonomous systems, in general, depend upon $t_0$.

## 3. Motivating Example

This section demonstrates application of our approach to the design of successive approximation to the norms of solutions to a Van der-Pol – like model with variable coefficients which was also tackled in [43] by our former technique. This system can be written as (1.1) and assumes the following form in dimensionless variables,

$$\frac{d}{dt} \begin{bmatrix} x_1 \\ x_2 \end{bmatrix} = \begin{pmatrix} 0 & 1 \\ -\omega^2 & -\alpha_1 \end{pmatrix} \begin{bmatrix} x_1 \\ x_2 \end{bmatrix} + \begin{bmatrix} 0 \\ -\alpha_2 x_2^3 \end{bmatrix} + \begin{bmatrix} 0 \\ F_2(t) \end{bmatrix} \tag{3.1}$$

where $x = [x_1 \; x_2]^T$, $f = [0 \; -\alpha_2 x_2^3]^T$, $F(t) = [0 \; F_2(t)]^T$, $\omega^2(t) = \omega_0^2 + \omega_1(t)$, $\omega_0 = const$, $\omega_1 = a_1 \sin r_1 t + a_2 \sin r_2 t$, $F_2 = a \sin \omega_2 t$, $a$, $\omega_2$, $a_i$, $r_i = const$, $i = 1, 2$. Next, we set in (1.11) that $Y_m = \|Y_{m1} \; Y_{m2}\|^T$, $y_i = [y_{1i} \; y_{2i}]^T$, $i \in \mathbb{N}_m$, where $\mathbb{N}_m$ is a set of $m$-positive integers, $y_1(t_0) = x_0$ and $y_i(t_0) = 0$, $i > 1$.

Then, substitution of (1.9) and (1.11) into (3.1) returns the following system,

$$\sum_{i=1}^{m} \dot{y}_{1i} + \dot{z}_1 = \sum_{i=1}^{m} y_{2i} + z_2$$

$$\sum_{i=1}^{m} \dot{y}_{2i} + \dot{z}_2 = -\omega_0 \left( \sum_{i=1}^{m} y_{1i} + z_1 \right) - \alpha_1 \left( \sum_{i=1}^{m} y_{2i} + z_2 \right) - \alpha_2 \left( \sum_{i=1}^{m} y_{2i} + z_2 \right)^3 + F(t)$$

Below we present two natural approaches to resolve this system which embrace the following conditions: equations for every $y_i$ are linear and coupled consecutively such that the prior sets of equations are independent from subsequent ones. This yields two sets of consecutive linear equations for $y_i$, $i \in \mathbb{N}_m$ and nonlinear equations for $z$ which are simulated in section 8.

Assume initially that $m = 1$ and $x = y_1 + z$. This comprises the following systems of linear equations for $y_1$ and nonlinear equations for $z$,

$$\dot{y}_{11} = y_{21}$$
$$\dot{y}_{21} = -\omega_0 y_{11} - \alpha_1 y_{21} + F(t), \quad y_1(t_0) = x_0 \tag{3.2}$$

$$\dot{z}_1 = z_2$$
$$\dot{z}_2 = -\omega_0 z_1 - \alpha_1 z_2 - \alpha_2 \left( 3 y_{21}^2 z_2 + 3 y_{21} z_2^2 + z_2^3 + y_{21}^3 \right), \quad z(t_0) = 0 \tag{3.3}$$

Then, application of equation (2.16) to (3.3) yields a scalar equation,

$$\dot{z}_{e1} = p(t) z_{e1} + c(t) |\alpha_2| \left( 3 |y_{21}^2| z_{e1} + 3 |y_{21}| z_{e1}^2 + z_{e1}^3 + y_{21}^3 \right), \quad z_{e1}(t_0) = 0 \tag{3.4}$$

where symbol $|.|$ stands for absolute value, both $p(t)$ and $c(t)$ are defined by (2.9) and (2.10) and $w(t)$ is a fundamental matrix of solutions to linearized and homogeneous equation (3.1). Clearly, $z_{e1} \geq \|z\|$.

For $m = 2$, $x = y_1 + y_2 + z$ and the corresponding system of equations includes (3.2) which is appendant by the following equations,

$$\dot{y}_{12} = y_{22}$$
$$\dot{y}_{22} = -\omega^2 y_{12} - \left( \alpha_1 + 3\alpha_2 y_{21}^2 \right) y_{22} - \alpha_2 y_{21}^3, \quad y_2(t_0) = 0 \tag{3.5}$$

$$\dot{z}_1 = z_2, \quad z(t_0) = 0,$$
$$\dot{z}_2 = -\omega_0 z_1 - \alpha_1 z_2 - \alpha_2 \left( 3 Y_{22}^2 z_2 + 3 Y_{22} z_2^2 + z_2^3 + (3 y_{21} + y_{22}) y_{22}^2 \right) \tag{3.6}$$

where we recall that $Y_{22} = y_{21} + y_{22}$.

Next, application of (2.16) to (3.6) returns the following equation,

$$\dot{z}_{e2} = p(t) z_{e2} + c(t) |\alpha_2| \left( 3 Y_{22}^2 z_{e2} + 3 |Y_{22}| z_{e2}^2 + z_{e2}^3 + |3 y_{21} + y_{22}| y_{22}^2 \right), \quad z_{e2}(t_0) = 0 \tag{3.7}$$

where $z_{e2} \geq \|z\|$.

Finally, for $m = 3$, the successive approximations include (3.2) and (3.5) which are appendant by the following equations,

$$\dot{y}_{13} = y_{23}, \qquad y_3(t_0) = 0,$$

$$\dot{y}_{23} = -\omega^2 y_{13} - \left(\alpha_1 + 3\alpha_2 Y_{22}^2\right) y_{23} - \alpha_2 \left(3 y_{21} + y_{22}\right) y_{22}^2 \qquad (3.8)$$

$$\dot{z}_1 = z_2, \quad z(t_0) = 0,$$

$$\dot{z}_2 = -\omega_0 z_1 - \alpha_1 z_2 - \alpha_2 \left(3 Y_{23}^2 z_2 + 3 Y_{23} z_2^2 + z_2^3 + (3 Y_{22} + y_{23}) y_{23}^2\right) \qquad (3.9)$$

where $Y_{23} = y_{21} + y_{22} + y_{23}$.

Consequently, the upper estimate of the norm of solutions to (3.9) assumes the following equation,

$$\dot{z}_{e3} = p(t) z_{e3} + c(t) |\alpha_2| \left(3 Y_{23}^2 z_{e3} + 3 |Y_{23}| z_{e3}^2 + z_{e3}^3 + |3 Y_{22} + y_{23}| y_{23}^2\right), \quad z_{e3}(t_0) = 0 \qquad (3.10)$$

Next, we outline an alternative approach where the underlying matrix of the linear block for all equations in successive approximations equals $A(t)$. To shorten the presentation, we present below only the equations for $m = 3$ that append (3.2) holding in this case,

$$\dot{y}_{12} = y_{22}$$

$$\dot{y}_{22} = -\omega^2 y_{12} - \alpha_1 y_{22} - \alpha_2 y_{21}^3, \quad y_2(t_0) = 0 \qquad (3.11)$$

$$\dot{y}_{13} = y_{23},$$

$$\dot{y}_{23} = -\omega^2 y_{13} - \alpha_1 y_{23} - \alpha_2 \left(Y_{22}^3 - y_{21}^3\right), \quad y_3(t_0) = 0 \qquad (3.12)$$

$$\dot{z}_1 = z_2, \quad z(t_0) = 0,$$

$$\dot{z}_2 = -\omega_0 z_1 - \alpha_1 z_2 - \alpha_2 \left(3 Y_{23}^2 z_2 + 3 Y_{23} z_2^2 + z_2^3 + Y_{23}^3 - Y_{22}^3\right) \qquad (3.13)$$

Afterwards, the upper estimate of the norm of solutions to (3.13) assumes the following equation,

$$\dot{z}_{e4} = p(t) z_{e4} + c(t) |\alpha_2| \left(3 Y_{23}^2 z_{e4} + 3 |Y_{23}| z_{e4}^2 + z_{e4}^3 + |Y_{23}^3 - Y_{22}^3|\right), \quad z_{e4}(t_0) = 0 \qquad (3.14)$$

Hence the transition matrices for equations (3.2), (3.11) and (3.12) are the same which simplifies their numerical evaluation. Nonetheless, simulations show that for a fixed $m$ the first set of equations is more accurate than the second one, see section 8.

Clearly, the above procedures can be readily applied to derivation of the underlying equations of our methodology for system (1) if $f(t, x)$ is a polynomial in $x$.

## 4. Successive Approximations

This section derives the successive approximations in a general form which does not assume that $f(t, x)$ is a polynomial in $x$. Yet, for polynomial systems, the derived equations reduce to ones that can be set up by extending the prior approach.

Adopting (1.9) and (1.11) for (1) let us to cast equations for $y_i(t, x_0)$ and $z(t, x_0)$ into two different forms. In the first case we assume that for sufficiently small $\|x\|$ and $\forall t \geq t_0$, Jacobean of $f$ in $x$, i.e., $f'(t, x) \in \mathbb{R}^{n \times n}$ is a continuous matrix in $t$ and $x$. Then, the first set of recursive equations for $y_k(t)$ can be written as follows,

$$\dot{y}_1 = Ay_1 + F(t), \ y_1(t_0) = x_0,$$
$$\dot{y}_2 = (A + f'(t, y_1))y_2 + f(t, y_1), \ y_2(t_0) = 0,$$
$$\dot{y}_3 = (A + f'(t, y_1 + y_2))y_3 + f(t, y_1 + y_2) - f(t, y_1) - f'(t, y_1)y_2, \ y_3(t_0) = 0 \quad (4.1)$$
$$\dot{y}_k = (A + f'(t, Y_{k-1}))y_k + f(t, Y_{k-1}) - f(t, Y_{k-2}) - f'(t, Y_{k-2})y_{k-1}, \ y_k(t_0) = 0,$$
$$3 \leq k \in \mathbb{N}_m$$

$$\dot{z} = Az + f(t, z + Y_m) - f(t, Y_{m-1}) - f'(t, Y_{m-1})y_m, \ z(t_0) = 0, \quad (4.2)$$
$$m \geq 2$$

Note that for $m = 1$, (4.2) takes the following form,
$$\dot{z} = Az + f(t, z + y_1), \ z(t_0) = 0$$

The second set of equations, generalizing our second approach from prior section, can be written as follows,
$$\dot{y}_1 = Ay_1 + F(t), \ y_1(t_0) = x_0,$$
$$\dot{y}_2 = Ay_2 + f(t, y_1), \ y_2(t_0) = 0, \quad (4.3)$$
$$\dot{y}_k = Ay_k + f(t, Y_{k-1}) - f(t, Y_{k-2}), \ y_k(t_0) = 0, \ 3 \leq k \in \mathbb{N}_m$$

$$\dot{z} = Az + f(t, z + Y_m) - f(t, Y_{m-1}), \ z(t_0) = 0, \ m \geq 2 \quad (4.4)$$

As in prior section, the underlying matrices in linear recursive equations (4.1) depend upon solutions of the preceding equations and, in turn, upon $t_0$ and $x_0$, $\forall k \geq 2$, whereas the underlying matrix in (4.3) equals $A$, $\forall k \geq 1$. Hence, for (4.3) the transition matrix is a block-diagonal which expedites repeated computations of solution to these equations emanated from a set of initial vectors.

Note that for $m = 1$ the corresponding equations are the same in both cases. Also in both sets of equations $y_k = y_k(t, t_0, x_0)$ and $z = z(t, t_0, x_0)$ since $y_k$, $k > 1$ and $z$ implicitly depend upon both $t_0$ and $x_0$ as parameters while complying with zero initial conditions. To shorten the notation we will frequently write below that $y_k = y_k(t, x_0)$ and $z = z(t, x_0)$.

## 5. Behavior of Approximate Solutions on Long Time-Intervals

This section examines behavior of $y_k(t, x_0)$, on long time-intervals playing a key role in our subsequent analysis. Let us assume firstly that $y_k(t, x_0)$ are defined by (4.1) and write the homogeneous counterpart of these linear equations as follows,
$$\dot{y}_k = (A - f'(t, Y_{k-1}))y_k, \ 1 \leq k \in \mathbb{N}_m \quad (5.1)$$

Next, we set that $\kappa_k = -v_k$ be the upper general exponents for (5.1) and $W_k(t, t_0)$, $k \in \mathbb{N}_m$ are the transition matrices for these equations. Hence, $W_1(t, t_0) = W(t, t_0)$, and $\|W_k(t, t_0)\| \leq N_k e^{-v_k(t-t_0)} \ \forall t \geq t_0, \ N_k > 0$.

### 5.1. Application of Lipschitz Continuity Condition

Let us assume in this subsection that the following Lipschitz - like continuity conditions hold for $t \geq t_0$,
$$\|f(t, s_1) - f(t, s_2)\| \leq l_2(t)\|s_1 - s_2\|, \ s_1, s_2 \in \Omega_3 \subset \mathbb{R}^n \quad (5.2)$$

$$\|f'(t,s)\| \leq l_3(t)\|s\|, \quad s \in \Omega_4 \subset \mathbb{R}^n \tag{5.3}$$

where $\Omega_i$, $i = 3,4$ are some bounded neighborhoods of $x \equiv 0$ and $l_k(t) \leq \hat{l}_k \in \mathbb{R}$, $t \geq t_0$, $k = 2,3$.

To simply further notations we set that $\Omega = \Omega_1 \cap \Omega_2 \cap \Omega_3 \cap \Omega_4$, where $\Omega$ is a bounded neighborhood of $x \equiv 0$. Next, we acknowledge some essential properties of $y_k(t,x_0)$ and exponents $v_k$ in the following theorems.

**Theorem 1**. Assume that $F_0 = 0$, $\|x_0\|$ is sufficiently small value, $v_1 > 0$, $y_k(t,x_0)$, $k \geq 1$ are solution to (4.1), inequalities (1.4), (5.2) and (5.3) hold, and $f'(t,x) \in \mathbb{R}^{n \times n}$ is a continuous matrix in $t$ and $x$. Then, $v_1 = v_k$, $\lim_{t_0 \leq t \to \infty} \|y_k(t,x_0)\| = 0$, $\|y_k(t,x_0)\| \leq \hat{y}_k(t,\|x_0\|) = O(\|x_0\|)$ and $\hat{y}_k$ increases in $\|x_0\|$ for $\forall t \geq t_0$.

**Proof**. We prove this statement by induction. Let $k = 1$, then, since $F_0 = 0$, $y_1(t,x_0) = W_1(t,t_0)x_0$ and $\|y_1(t,x_0)\| \leq \hat{y}_1(t,\|x_0\|) = \|x_0\|N_1 \exp(-v_1(t-t_0))$, $\forall t \geq t_0$. Hence, $\lim_{t \to \infty} \|y_1(t,x_0)\| = 0$, $y_1(t,x_0) \in \Omega$, $\forall t \geq t_0$ for sufficiently small $\|x_0\|$ and, due to (5.3), $\|f'(t,y_1(t,x_0))\| \leq \hat{l}_3 \|y_1(t,x_0)\|$, which implies that $\lim_{t_0 \leq t \to \infty} \|f'(t,y_1(t,x_0))\| = 0$. Thus, due to Statement 2, $v_1 = v_2$.

Let us next bound $y_2(t,x_0)$ and $y_3(t,x_0)$ to unveil common traits in this sequence. In fact, $\|W_2(t,t_0)\| \leq N_2 e^{-v_1(t-t_0)}$, $\|f(t,y_1)\| \leq \hat{l}_1 \|y_1\|$ for sufficiently small $\|x_0\|$, and

$$\|y_2(t,x_0)\| \leq \int_{t_0}^{t} \|W_2(t,\tau)\| \|f(t,y_1(\tau,x_0))\| d\tau \leq \hat{l}_1 \int_{t_0}^{t} \|W_2(t,\tau)\| \|y_1(\tau,x_0)\| d\tau.$$ Substitution of the prior bounds in the last integral yields that,

$\|y_2(t,x_0)\| \leq \hat{y}_2(t,\|x_0\|) = \mu_{21} \|x_0\| e^{-v_1(t-t_0)}(t-t_0)$, $\forall t \geq t_0$, $0 < \mu_{21} = const$, which implies that $\lim_{t_0 \leq t \to \infty} \|y_2(t,x_0)\| = 0$.

Next, we infer that, $y_1(t,x_0) + y_2(t,x_0) \in \Omega$, $\forall t \geq t_0$ for sufficiently small $\|x_0\|$ and, due to (5.3), $\|f'(t,y_1(t,x_0)+y_2(t,x_0))\| \leq \hat{l}_3(\|y_1(t,x_0)\| + \|y_2(t,x_0)\|)$. The last inequality implies that, $\lim_{t_0 \leq t \to \infty} \|f'(t,y_1(t,x_0)+y_2(t,x_0))\| = 0$ which with Statement 2 yields that $v_3 = v_1$.

In turn, application of (5.2) and (5.3) infers that $\|f(y_1+y_2) - f(y_1)\| \leq \hat{l}_1 \|y_2\|$ and $\|f'(y_1)y_2\| \leq \hat{l}_1 \hat{l}_2 \|y_1\| \|y_2\|$. Then, application of standard integration formulas returns that,

$\|y_3(t,x_0)\| = \hat{y}_3(t,\|x_0\|) = \mu_{31} \|x_0\| e^{-v_1(t-t_0)}(t-t_0)^2 + \|x_0\|^2 \left(\mu_{32} e^{-2v_1(t-t_0)}(t-t_0) + \mu_{33} e^{-v_1(t-t_0)}\right)$, $\forall t \geq t_0$, $0 < \mu_{3i} = const$, $i \in \mathbb{N}_3$. Hence, $\lim_{t_0 \leq t \to \infty} \|y_3(t,x_0)\| = 0$.

Afterward, to simplify the induction let us define functions, $\alpha_{sj}(t-t_0) = e^{-sv_1(t-t_0)} p_j(t-t_0)$, $\forall t \geq t_0$, where $s \in \mathbb{N}$, $j \in \mathbb{Z} = 0 \cup \mathbb{N}$ and $p_j(t-t_0)$ is a polynomial in $t-t_0$ of degree $j$ and $p_0 = const$ which can be zero. Clearly, $\lim_{t_0 \leq t \to \infty} \alpha_{sj}(t-t_0) = 0$.

Next, let $\Theta_1$ is a vector space, over $\mathbb{R}$ which is spanned by $\alpha_{sj}(t-t_0)$ and also closed with respect to application of two operations: multiplication of its entries, i.e., $\alpha_{s_1 j_1} \alpha_{s_2 j_2} \in \Theta_1$, $s_i \in \mathbb{N}$, $j_i \in \mathbb{Z}$, $i = 1, 2$ and a convolution, which is defined as follows, $\int_{t_0}^{t} e^{-v_1(t-\tau)} \alpha_{sj}(\tau - t_0) d\tau \in \Theta_1$. Note that the last formula follows from the standard integration rules.

In turn, we assume by induction that,

$$\|y_k(t, x_0)\| \leq \hat{y}_k(t, \|x_0\|) = \sum_{i=1}^{F_k} \|x_0\|^i \rho_{ki}(t - t_0), \forall t \geq t_0, \forall k \in \mathbb{N}_K \quad (5.4)$$

where $F_k$ is $k$-th Fibonacci number, $0 < \rho_{ki} \in \Theta_1$, $\forall t > t_0$. Clearly, $\hat{y}_k(t, \|x_0\|) \in \Theta_1$ and increases in $\|x_0\|$, $\forall t \geq t_0$ which implies that $\lim_{t_0 \leq t \to \infty} \|y_k(t, x_0)\| = 0$.

Next, we show in sequence that (5.4) holds for $k = K+1$. In fact, repetition of the foregoing steps yields that, $\lim_{t \to \infty} \|f'(t, Y_K)\| = 0$, which implies that, $v_{K+1} = v_1$ and $\|W_{K+1}(t, t_0)\| \leq N_{K+1} e^{-v_1(t-t_0)}$. Next, as prior, we infer that $Y_K \in \Omega$ for sufficiently small $\|x_0\|$ and (5.2) and (5.3) hold under this condition. This implies that,

$\|f(t, Y_K) - f(t, Y_{K-1})\| \leq \hat{l}_2 \|y_K\| \leq \hat{l}_2 \hat{y}_K \in \Theta_1$ and $\|f'(t, Y_{K-1}) y_K\| \leq \hat{l}_3 \|y_K\| \|Y_{K-1}\| \leq l_3 \hat{y}_K \hat{Y}_{K-1} \in \Theta_1$, where $\hat{Y}_{K-1} = \sum_{j=1}^{K-1} \hat{y}_i$. Consequently, we get that, $\|y_{K+1}(t, x_0)\| \leq \hat{y}_{K+1}(t, \|x_0\|) =$

$N_{K+1} \int_{t_0}^{t} e^{-v_1(t-\tau)} \hat{y}_K(\tau, \|x_0\|)(\hat{l}_2 + \hat{l}_3 \hat{Y}_{K-1}(\tau, \|x_0\|)) d\tau \in \Theta_1, \forall t \geq t_0$. This implies that $\lim_{t_0 \leq t \to \infty} \|y_{K+1}(t, x_0)\| = 0$

and $\hat{y}_{K+1}(t, \|x_0\|)$ increases in $\|x_0\|$ since integral of positive functions, i.e., $\int_{t_0}^{t} e^{-v_1(t-\tau)} \rho_{ki}(\tau - t_0) d\tau > 0$ and

$\int_{t_0}^{t} e^{-v_1(t-\tau)} \rho_{ki}(\tau - t_0) \rho_{sj}(\tau - t_0) d\tau > 0$ □

**Theorem 2**. Assume that $F_0 > 0$, $v_1 > 0$ and both, $\|x_0\|$ and $F_0$ are sufficiently small, $y_k(t, x_0)$ is a solution to (4.1) and inequalities (1.4), (5.2) and (5.3) are intact and $f'(t, x) \in \mathbb{R}^{n \times n}$ is a continuous matrix in both $t$ and $x$. Then, $0 < v_k < v_{k-1}$, $k > 1$, $\lim_{t_0 \leq t \to \infty} \|y_k(t, x_0)\| \leq O(F_0)$, $\|y_k(t, x_0)\| \leq \bar{y}_k(t, \|x_0\|, F_0) = O(\|x_0\| + F_0)$, $t \geq t_0$ and $\bar{y}_k(t, \|x_0\|, F_0)$ increases in both $\|x_0\|$ and $F_0$, $\forall t \geq t_0$.

**Proof**. We proof this theorem by induction. In fact, $\|y_1(t, x_0)\| \leq \|x_0\| N_1 \exp(-v_1(t-t_0))$
$+ (N_1 F_0 / v_1)(1 - e^{-v_1(t-t_0)}) = O(\|x_0\| + F_0)$, $\forall t \geq t_0$ and $\lim_{t_0 \leq t \to \infty} \|y_1(t, x_0)\| \leq O(F_0)$. Thus, $y_1 \in \Omega$ and $\|f'(t, y_1)\| \leq \hat{l}_3 \|y_1\|$ if both $\|x_0\|$ and $F_0$ are sufficiently small. Hence, due to Statement 1, we infer that under appropriate choice of $\|x_0\|$ and $F_0$, $0 < v_2 < v_1$ which, subsequently, yields that,

$\|W_2(t,t_0)\| \leq N_2 e^{-v_2(t-t_0)}$, $\forall t \geq t_0$. Furthermore, since under our assumptions, $\|f(y_1)\| \leq \hat{l}_1 \|y_1\|$ we get that,

$$\|y_2\| \leq \int_{t_0}^{t} \|W_2(t,\tau)\| \|f(y_1(\tau,x_0))\| d\tau \leq \hat{l}_1 \int_{t_0}^{t} \|W_2(t,\tau)\| \|y_1(\tau,x_0)\| d\tau.$$

Next, we introduce some notations, which will be used afterward. Let us write the previous formula as follows,

$\|y_2(t,x_0)\| \leq \bar{y}_2(t,\|x_0\|,F_0) = \beta_{210}(t-t_0)\|x_0\| + \beta_{211}(t-t_0)F_0$, where

$$\beta_{210}(t-t_0) = N_1 N_2 \int_{t_0}^{t} e^{-v_2(t-\tau)} e^{-v_1(\tau-t_0)} d\tau, \; \beta_{211}(t-t_0) = N_1 N_2 v_1^{-1} \int_{t_0}^{t} e^{-v_2(t-\tau)} \left(1 - e^{-v_1(\tau-t_0)}\right) d\tau.$$ Clearly,

$\lim_{t_0 \leq t \to \infty} \beta_{210}(t-t_0) = 0$, $\lim_{t_0 \leq t \to \infty} \beta_{211}(t-t_0) = const \neq 0$, and $\beta_{2ij}(t-t_0) > 0$, $i = 1, j = 0,1$ since they are defined through integrations of the products of positive functions. Hence, $\lim_{t_0 \leq t \to \infty} \|y_2(t,x_0)\| \leq O(F_0)$,

$\bar{y}_2(t,\|x_0\|,F_0) = O(\|x_0\| + F_0)$, $t \geq t_0$ and the last function increases in both $\|x_0\|$ and $F_0$, $\forall t \geq t_0$.

Next, let $\Theta_2$ be a vector space over $\mathbb{R}$ that is spanned by functions, $\alpha_{sj}^*(t-t_0) = e^{-s(t-t_0)} p_j(t-t_0)$, $\forall t \geq t_0$, where $s \in \mathbb{R}_{>0}$, $j \in \mathbb{Z}$ and $p_j(t-t_0)$ is a polynomial in $t - t_0$ of degree $j$ with $p_0 = const$ which can be zero. Note that $\Theta_2$ is closed with respect to multiplication of its entries and a convolution which is defined as follows,

$$\int_{t_0}^{t} e^{-v_k(t-\tau)} \alpha_{sj}^*(\tau-t_0) d\tau, \; v_k > 0, \; k \in \mathbb{N}_m.$$

Next, we assume that,

$$\|y_k(t,x_0)\| \leq \bar{y}_k(t,\|x_0\|,F_0) = \sum_{i=1}^{F_k} \sum_{j=0}^{i} \beta_{kij}(t-t_0) \|x_0\|^{i-j} F_0^j, \; \forall k \in \mathbb{N}_K \tag{5.5}$$

where $\beta_{kij}(t-t_0) > 0$, $\beta_{kij}(t-t_0) \in \Theta_2$, $\lim_{t_0 \leq t \to \infty} \beta_{kij}(t-t_0) = 0$, $\forall i \neq j$, $\lim_{t_0 \leq t \to \infty} \beta_{kii}(t-t_0) = const \neq 0$, $\forall t \geq t_0$, $F_k$ is $k$-th Fibonacci number. Clearly, $\bar{y}_k(t,\|x_0\|,F_0) \in \Theta_2$, $\bar{y}_k(t,\|x_0\|,F_0) = O(\|x_0\| + F_0)$, $\forall t \geq t_0$, $\lim_{t_0 \leq t \to \infty} \|y_k(t,x_0)\| \leq O(F_0)$ and our assumption holds for $k = 2$. Let us show by induction that (5.5) holds for $k = K+1$. In fact, $y_k(t,x_0) \in \Omega$, $k \in \mathbb{N}_K$, $Y_K \in \Omega$ since both $\|x_0\|$ and $F_0$ are sufficiently small and

$\|Y_K\| \leq \sum_{k=1}^{K} \|y_k\| \leq \bar{Y}_K = \sum_{i=1}^{K} \bar{y}_i = O(\|x_0\| + F_0)$, $\forall t \geq t_0$. This implies that, $\|f'(t,Y_K)\| \leq \hat{l}_3 \bar{Y}_K \in \Theta_2$ and $\lim_{t_0 \leq t \to \infty} \|f'(t,Y_K)\| \leq O(F_0)$, which, with Statement 1, infers that $0 < v_{K+1} < v_K$.

Next, for sufficiently small $\|x_0\|$ and $F_0$, application of (5.2) and (5.3) implies that,

$\|f(t,Y_K) - f(t,Y_{K-1})\| \leq \hat{l}_2 \bar{y}_{K-1} \in \Theta_2$ and $\|f'(t,Y_{K-1}) y_K\| \leq \hat{l}_3 \bar{y}_K \bar{Y}_{K-1} \in \Theta_2$. In turn,

$$\|y_{K+1}(t,x_0)\| \leq \bar{y}_{K+1}(t,\|x_0\|,F_0) =$$

$$N_{K+1} \int_{t_0}^{t} \exp(-v_{K+1}(t-\tau)) \bar{y}_K(\tau,\|x_0\|,F_0) \left[\hat{l}_2 + \hat{l}_3 \bar{Y}_{K-1}(\tau,\|x_0\|,F_0)\right] d\tau$$

Consequently, $\|y_{K+1}(t,x_0)\| \leq \bar{y}_{K+1}(t,\|x_0\|,F_0) = \sum_{i=1}^{F_{K+1}}\sum_{j=0}^{i}\beta_{K+1,ij}(t-t_0)\|x_0\|^{i-j}F_0^j \in \Theta_2$, where $\beta_{K+1,ij}(t-t_0)$ are defined by prior formula, and $\bar{y}_{K+1} = O(\|x_0\|+F_0)$. Note that $\beta_{K+1,ij}(t-t_0) > 0, \forall t \geq t_0$ since multiplication and integration of positive functions yield positive functions and $\beta_{K+1,ij}(t-t_0) \in \Theta_2$ since $\Theta_2$ is closed with respect to both multiplication of its elements and prior defined convolution. Furthermore, application of standard integration rules yields that, $\lim_{t_0 \leq t \to \infty}\beta_{K+1,ij}(t-t_0) = 0, \forall i \neq j$ and $\lim_{t \to \infty}\beta_{K+1,ii}(t-t_0) = const \neq 0$.

Hence, $\bar{y}_{K+1}(t,\|x_0\|,F_0)$ increases in both $\|x_0\|$ and $\hat{F}$ and $\lim_{t_0 \leq t \to \infty}\|y_{K+1}(t,x_0)\| \leq O(F_0)$ □

Subsequently, we derive the upper bound for the norm of $y_k(t,x_0)$ which is determined by formulas (4.3) in the following,

**Theorem 3**. Assume that both $\|x_0\|$ and $F_0$ are sufficiently small, $v_1 > 0$ and $y_k(t,x_0)$, $k \geq 1$ are solution to (4.3) and inequalities (1.4), (5.2) and (5.3) hold. Then,

$$\|y_k(t,x_0)\| \leq \tilde{y}_k(t,\|x_0\|,F_0) =$$
$$\hat{l}^{k-1}N_1^k\left\{\left[(\|x_0\|-F_0/v_1)\frac{(t-t_0)^{k-1}}{(k-1)!} - ... - \frac{F_0}{v_1^k}\right]\exp(-v_1(t-t_0)) + \frac{F_0}{v_1^k}\right\}, k \geq 1, \forall t \geq t_0 \quad (5.6)$$

where $l = \max(\hat{l}_1,\hat{l}_2)$ and $\tilde{y}_k(t,\|x_0\|,F_0)$ increases in both $\|x_0\|$ and $F_0$, $\forall t \geq t_0$, and series $\sum_{k=1}^{\infty}y_k(t,x_0)$ converge absolutely $\forall t \geq t_0$ if $lN_1F_0/v_1 < 1$.

**Proof.** The first part of this statement can be readily proved by induction. In fact,

$$\|y_1(t,x_0)\| \leq \tilde{y}_1(t,\|x_0\|,F_0) = N_1\|x_0\|e^{-v_1(t-t_0)} + (N_1F_0/v_1)(1-e^{-v_1(t-t_0)}) =$$
$$= N_1((\|x_0\|-F_0/v_1)e^{-v_1(t-t_0)} + F_0/v_1)$$

Hence, for sufficiently small $\|x_0\|$ and $F_0$, $y_1(t,x_0) \in \Omega$, $\forall t \geq t_0$, which implies that under these conditions, $\|f(y_1)\| \leq \hat{l}\|y_1\|$.

In turn,

$$\|y_2(t,x_0)\| \leq \tilde{y}_2(t,\|x_0\|,F_0) = N_1\hat{l}\int_{t_0}^{t}e^{-v_1(t-\tau)}\tilde{y}_1(\tau,\|x_0\|,F_0)d\tau =$$
$$N_1^2\hat{l}\left[((\|x_0\|-F_0/v_1)(t-t_0)-(F_0/v_1^2))e^{-v_1(t-t_0)} + F_0/v_1^2\right]$$

and, subsequently, under the assumption of this statement, $y_2(t,x_0) \in \Omega$, $\forall t \geq t_0$. In sequence, we assume that $\|y_k(t,x_0)\|$, $2 < k \in \mathbb{N}_K$ are bounded from above by (5.6) and show that this formula holds for $k = K+1$. Clearly, in this case, both $y_K$ and $Y_K \in \Omega$ for sufficiently small $\|x_0\|$ and $F_0$. In turn,

$\|f(t,Y_K) - f(t,Y_{K-1})\| \leq \hat{l}\|y_K\| \leq \hat{l}\tilde{y}_K(t,\|x_0\|,F_0)$ and consequently, $\|y_{K+1}\| \leq \tilde{y}_{K+1}(t,\|x_0\|,F_0) =$

$N_1\hat{l}\int_{t_0}^{t}e^{-v_1(t-\tau)}\tilde{y}_K(\tau,\|x_0\|,F_0)d\tau$. Then, application of standard formulas infers the first part of this statement.

Next, let us show by induction that $\tilde{y}_k(t, \|x_0\|, F_0)$, $k \geq 1$ increases in both $\|x_0\|$ and $F_0$, $\forall t \geq t_0$. In fact, the statement holds for $k = 1$. Assume, subsequently, that $\tilde{y}_k(t, \|x_0\|, F_0) = \|x_0\|\eta_{k1}(t) + F_0\eta_{k2}(t)$, $2 < k \in \mathbb{N}_K$ where continuous $\eta_{ki}(t) > 0$, $t \geq t_0$, $i = 1, 2$. Then,

$$\tilde{y}_{K+1}(t, \|x_0\|, F_0) = N_1\hat{l}\int_{t_0}^{t} e^{-v_1(t-\tau)}\tilde{y}_K(\tau, \|x_0\|, F_0)d\tau =$$

$$N_1\hat{l}\left(\|x_0\|\int_{t_0}^{t} e^{-v_1(t-\tau)}\eta_{K1}(\tau)d\tau + F_0\int_{t_0}^{t} e^{-v_1(t-\tau)}\eta_{K1}(\tau)d\tau\right)$$

Obviously, the integrals of the products of positive functions entering the above formula are positive, which embraces the second part of this statement.

Finally, to prove absolute convergence of series, $\sum_{k=1}^{\infty} y_k(t, x_0)$, we apply the Cauchy's root test to (5.6). Clearly,

$\lim_{k\to\infty}\left[(t-t_0)^{k-1}/(k-1)!\right] = 0$ and $\exp(-v_1(t-t_0)) \leq 1$, $\forall t \geq t_0$ which simplifies the test and infers our statement □

### 5.2 Application to Polynomial Systems

In appendix of this section we adopted for simplicity the notation used prior in subsection 5.1 in the matching statements. Note that functions, $\hat{y}_k(t, \|x_0\|)$, $\overline{y}_k(t, \|x_0\|, F_0)$ and $\tilde{y}_k(t, \|x_0\|, F_0)$ play the same role in both subsections but their definitions in this subsection are different from the prior ones.

Let us assume next that $f(t, x)$ with $f(t, 0) = 0$ is a vector – polynomial in $x$ without linear components. Then application of inequalities (1.4), (5.2) and (5.3) to theorems 1 -3 turns out to be redundant which let to void a previously essential constrain that $\|x_0\|$ is sufficiently small. This leads to the following counterparts of theorems 1 - 3.

**Theorem 4.** Assume that $f(t, x)$ with $f(t, 0) = 0$ is a vector – polynomial in $x$ without linear components, $F_0 = 0$, $v_1 > 0$ and $y_k(t, x_0)$, $k \geq 1$ are solution to (4.1). Then, $v_1 = v_k$, $\lim_{t\to\infty}\|y_k(t, x_0)\| = 0$, $\|y_k(t, x_0)\| \leq \hat{y}_k(t, \|x_0\|) = O(\|x_0\|)$, $t \geq t_0$ and $\hat{y}_k$ increases in $\|x_0\|$, $\forall t \geq t_0$, $k \geq 1$.

**Proof.** The proof of this statement is analogous to the previous one with some exceptions. The initial step of induction and formula (5.4) are intact in this case but $F_k$ here is the number of additions in (5.4). To proof the induction, we need to show that under the assumptions of this theorem, $\lim_{t_0 \leq t \to \infty}\|f'(t, Y_K)\| = 0$ and

$\|f(t, Y_K)\| \leq f_*(t, \hat{Y}_K) \in \Theta_1$, where $\|Y_K\| \leq \hat{Y}_K = \sum_{i=1}^{K}\hat{y}_i(t, \|x_0\|) \in \Theta_1$, $f_*(t, \hat{Y}_K)$ is a scalar polynomial in $\hat{Y}_K$

and $f_*(t, 0) = 0$. In fact, these inferences directly follow from repeated applications of standard inequalities that are briefly discussed in section 2.3, see also Appendix of this paper for the similar derivations.

Firstly, we note that the entries of Jacobian matrix, $f'(t, Y_K)$, i.e. $f_{1,ij}(t, Y_K)$, $i, j \in \mathbb{N}_n$ are polynomials in $Y_K$ of degrees greater or equal to one with $f_{1,ij}(t, 0) = 0$. Using mentioned above inequalities we infer that

$\|f_{1ij}(t, Y_K)\| \leq f_{2,ij}(t, \|Y_K\|)$, where $f_{2,ij}$ are scalar polynomials in $\|Y_K\|$ of degrees greater or equal to one and $f_{2,ij}(t,0) = 0$. In turn, $\|f'(t,Y_K)\| \leq \|f'(t,Y_K)\|_F \leq \left(\sum_{i,j=1}^{n} f_{2,ij}^2(t,\|Y_k\|)\right)^{1/2}$ and, due to assumption of the induction, $\lim_{t_0 \leq t \to \infty} \|Y_K\| \leq \lim_{t \to \infty} \hat{Y}_K = \lim_{t_0 \leq t \to \infty} \sum_{i=1}^{K} \hat{y}_i(t,\|x_0\|) = 0$, which implies that $\lim_{t_0 \leq t \to \infty} \|f'(t,Y_K)\| = 0$.

Next, application of the mentioned above inequalities yields that, $\|f(t,Y_K)\| \leq \|f(t,Y_K)\|_1 = \sum_{i=1}^{n} |f_i(t,Y_K)|$
$\leq f_*(t,\|Y_K\|) \leq f_*(t,\hat{Y}_K)$, where $f_*$ is a scalar polynomial in $\hat{Y}_K$ with degree greater or equal two. Hence, $f_*(t,\hat{Y}_K) \in \Theta_1$.

Obviously, $f_*(t,\hat{Y}_{K-1}) \in \Theta_1$ as well and $\|f'(t,Y_{K-1})y_K\| \leq \hat{y}_K f_3(t,\hat{Y}_{K-1}) \in \Theta_1$, where $f_3$ is a scalar polynomial in $\hat{Y}_{K-1}$ □

**Theorem 5.** Assume that $f(t,x)$ is a vector – polynomial in $x$ without linear components and $f(t,0) = 0$, $F_0 > 0$ is a sufficiently small value, $v_1 > 0$ and $y_k(t,x_0)$, $k \geq 1$ are solution to (4.1). Then,
$0 < v_k < v_{k-1}$, $k > 1$, $\lim_{t_0 \leq t \to \infty} \|y_k(t,x_0)\| = O(F_0)$, $\|y_k(t,x_0)\| \leq \bar{y}_k(t,\|x_0\|,F_0) = O(\|x_0\| + F_0)$, $t \geq t_0$
and $\bar{y}_k(t,\|x_0\|,F_0)$ increases in both $\|x_0\|$ and $F_0$, $\forall t \geq t_0$, $k \geq 1$.

**Proof.** The proof of this statement is analogous to proofs of theorems 2 and 4. In fact, the initial step of induction and formula (5.5) are intact but, as prior, in this case $F_k$ is the number of terms in (5.5). To embrace the induction, we need to show that under the assumptions of this theorem, $\lim_{t_0 \leq t \to \infty} \|f'(t,Y_k)\| = O(F_0)$ and

$\|f(t,Y_k)\| \leq f^*(t,\bar{Y}_k) \in \Theta_2$, where $\|Y_k\| \leq \bar{Y}_k = \sum_{i=1}^{k} \bar{y}_i(t,\|x_0\|,F_0) \in \Theta_2$, $f^*(t,\bar{Y}_k)$ is a scalar polynomial in $\bar{Y}_k$, $f^*(t,\bar{Y}_k) = O(\|x_0\| + F_0)$, $\forall t \geq t_0$ and $\lim_{t \to \infty} f^*(t,\bar{Y}_k) = O(F_0)$. The proofs of these statements are similar to ones made in theorem 4. In fact, the induction assumes that, $\lim_{t_0 \leq t \to \infty} \bar{y}_k(t,\|x_0\|,F_0) = O(F_0)$, $\bar{y}_k(t,\|x_0\|,F_0) = O(\|x_0\| + F_0)$, $\forall t \geq t_0$, $k \in \mathbb{N}_K$ and $\bar{y}_k(t,\|x_0\|,F_0)$ increases in $\|x_0\|$ and $F_0$ $\forall t \geq t_0$, $k \in \mathbb{N}_K$. Thus, application of the steps used to prove the previous statements infers this theorem □

**Theorem 6.** Assume that $f(t,x)$ with $f(t,0) = 0$ is a vector – polynomial in $x$, $v_1 > 0$ and $y_k(t,x_0)$, $k \geq 1$ are solution to (4.3). In addition, we assume that: 1. $F_0 = 0$ or 2. $F_0 > 0$. Then, 1. $\lim_{t \to \infty} \|y_k(t,x_0)\| = 0$, $\|y_k(t,x_0)\| \leq \hat{y}_k(t,\|x_0\|) = O(\|x_0\|)$ and $\hat{y}_k$ increases in $\|x_0\|$ for $\forall t \geq t_0$. 2. $\lim_{t \to \infty} \|y_k(t,x_0)\| \leq O(F_0)$, $\|y_k(t,x_0)\| \leq \tilde{y}(t,\|x_0\|,F_0) = O(\|x_0\| + F_0)$ and $\tilde{y}(t,\|x_0\|,F_0)$ increases in both $\|x_0\|$ and $F_0$, $\forall t \geq t_0$.

**Proof.** In fact, in this case the 'linear block' of (4.3) is defined by matrix $A(t)$ which immediately yields that $v_k = v_1$, $k > 1$. Thus, we merely need to show that: 1. $\|f(t,Y_K)\| \leq f_+(t,\hat{Y}_K) \in \Theta_1$ if $F_0 = 0$ and 2. $\|f(t,Y_K)\| \leq f^+(t,\hat{Y}_K) \in \Theta_2$ if $F_0 > 0$, where both $f_+$ and $f^+$ are scalar polynomials in $\hat{Y}_K$. The proofs of these statements are directly followed from mentioned above inequalities and are completely analogous to the prior ones □

Clearly, under conditions of theorems 1-6, $\lim_{t_0 \le t \to \infty} \|x(t, x_0)\| = \lim_{t_0 \le t \to \infty} \|z(t, x_0)\|$ if $F_0 = 0$ and

$\lim_{t_0 \le t \to \infty} \|x(t, x_0)\| = \lim_{t_0 \le t \to \infty} \|z(t)\| + O(F_0)$ if $F_0 > 0$. Hence, behavior of the norm of error term, $\|z(t)\|$ on long time intervals is circumstantial in analysis of boundedness and stability of solutions to (1) and (2), respectively.

## 6. Linear Error Equation and Its Applications

This section applies the methodology developed in [43] to estimate the error norm of approximate solutions given by (4.1), i.e., the norm of solution to (4.2) using Lipschitz continuity condition. This, in turn, yields successive approximations to bilateral solution bounds and boundedness/stability conditions for solutions of equations (1) and (2). Our inferences can be readily attuned if the error term is defined by equation (4.4) as well.

### 6.1. Linear Error Equation

Firstly, we apply inequality (2.8) to equation (4.2) as follows,

$$D^+\|z_1\| \le p(t)\|z_1\| + c(t)\|f(t, z_1 + Y_m) - f(t, Y_{m-1}) - f'(t, Y_{m-1})y_m\|, \ \forall t \ge t_0, \ \|z_1(t_0)\| = 0 \qquad (6.1)$$

where functions, $p(t)$ and $c(t) \le \hat{c} < \infty$ are defined in section 2, see also [43]. Clearly, $\|z_1\| \ge \|z\|$ and both, $Y_m \in \Omega$ and $Y_{m-1} \in \Omega$, $\forall t \ge t_0$ if $\|x_0\|$ and $F_0$ are sufficiently small. We will show below that this condition also yields that, $z_1 + Y_m \in \Omega$, $\forall t \ge t_0$ as well. Hence, application of standard norm's inequality and (5.2) transforms (6.1) into the following form,

$$D^+\|z_1\| \le \lambda(t)\|z_1\| + \Gamma(t, x_0), \ \forall t \ge t_0$$
$$\|z_1\|(t_0, x_0) = 0$$

where

$$\Gamma(t, x_0) = c(t)\left(l_2(t)\|y_m(t, x_0)\| + \|f'(t, Y_{m-1})y_m(t, x_0)\|\right) \qquad (6.2)$$

and $\lambda(t) = p(t) + c(t)l_2(t)$. Subsequent application of comparison principle to the last inequality leads to the following linear equation,

(6.3)
$$\dot{Z}_1 = \lambda(t)Z_1 + \Gamma, \ \forall t \ge t_0$$
$$Z_1(t_0, x_0) = 0$$

where, $Z_1(t, x_0) \ge \|z_1(t, x_0)\|$, $t \ge t_0$. Let us also recall that both $z_1(t, x_0)$ and $Z_1(t, x_0)$ depend implicitly upon both $t_0$ and $x_0$ but assume zero initial condition at $t = t_0$.

Due to our initial assumptions, $p(t)$, $c(t)$, $l_2(t)$ and, hence, $\lambda(t)$ are continuous and bounded functions in $t$. In addition, under assumptions of either Theorems 1 ($F_0 = 0$) or Theorem 2 ($F_0 > 0$), we infer that continuous $\Gamma \le \infty$, $\forall t \ge t_0$. In the latter case we infer that,

$$\Gamma \le \overline{\Gamma}(t, \|x_0\|, F_0) = \hat{c}\overline{y}_m(t, \|x_0\|, F_0)(\hat{l}_2 + \hat{l}_3\overline{Y}_{m-1}(t, \|x_0\|, F_0)) \qquad (6.4)$$

and under assumptions of Theorem 2, we get that, $\overline{\Gamma}(t, \|x_0\|, F_0) \in \Theta_2$, $\lim_{t_0 \le t \to \infty} \overline{\Gamma}(t, \|x_0\|, F_0) = O(F_0)$, $\overline{\Gamma}(t, \|x_0\|, F_0)$ increases in both $\|x_0\|$ and $F_0$, $\forall t \ge t_0$ and, $\sup_{t_0 \le t} \overline{\Gamma}(t, \|x_0\|, F_0) = O(\|x_0\| + F_0)$ since $\sup_{t_0 \le t} \overline{y}(t, \|x_0\|, F_0) = O(\|x_0\| + F_0)$.

In turn, for $F_0 = 0$ (6.4) should be attuned as follows,

$$\Gamma \leq \hat{\Gamma}(t, F_0, \|x_0\|) = \hat{k}\hat{y}_m(t, \|x_0\|)(\hat{l}_2 + \hat{l}_3 \hat{Y}_{m-1}(t, \|x_0\|, F_0)) \tag{6.5}$$

and under the assumptions of Theorem 1, $\hat{\Gamma}(t, \|x_0\|) \in \Theta_1$, $\lim_{t_0 \leq t \to \infty} \hat{\Gamma}(t, \|x_0\|) = 0$, $\sup_{t_0 \leq t} \hat{\Gamma}(t, x_0) = O(\|x_0\|)$ and $\hat{\Gamma}(t, x_0)$ increases in $\|x_0\|$, $\forall t \geq t_0$.

Subsequently, the solution to (6.3) can be written as,

$$Z_1(t, x_0) = \int_{t_0}^{t} \exp\left(\int_{\tau}^{t} \lambda(s, t_0) ds\right) \Gamma(\tau, x_0) d\tau, \forall t \geq t_0 \tag{6.6}$$

### 6.2. Application of Linear Error Equation

Thereafter, we present the following theorems which provide bilateral bounds and corresponding boundedness/stability criteria for solutions to (1) or (2), respectively.

**Theorem 7**. Assume that booth, $\|x_0\|$, $F_0 > 0$ are sufficiently small, $\lambda(t) < -\hat{\lambda}$, $\forall t \geq t_0$, $\hat{\lambda} > 0$, $f'(t, x)$ is a continuous matric, inequalities (1.4), (5.2) and (5.3) hold and $c(t) \leq \hat{c}$. Then, $\|x(t, x_0)\| < \infty$, $\forall t \geq t_0$, $\lim_{t_0 \leq t \to \infty} \|x(t, x_0)\| = O(F_0)$, where $x(t, x_0)$ is a solution to (1) and

$$\|Y_m(t, x_0)\| - \bar{Z}_1(t, \|x_0\|, F_0) \leq \|Y_m(t, x_0)\| - Z_1(t, x_0) \leq \|x(t, x_0)\| \leq \|Y_m(t, x_0)\| + Z_1(t, x_0)$$
$$\leq \bar{Y}_m(t, \|x_0\|, F_0) + \bar{Z}_1(t, \|x_0\|, F_0), \forall t \geq t_0 \tag{6.7}$$

where $\bar{Z}_1(t, \|x_0\|, F_0) = \int_{t_0}^{t} e^{-\hat{\lambda}(t-\tau)} \bar{\Gamma}(\tau, \|x_0\|, F_0) d\tau$ and both left – side inequalities should be adjusted to zero if they assume negative values.

**Proof**. In fact, at the end of section 2, we indicated that, $v_1 \geq \hat{\lambda} > 0$. Then, due to Theorem 2,

$$\|Y_m(t, x_0)\| \leq \bar{Y}_m = \sum_{k=1}^{m} \bar{y}_k(t, \|x_0\|, F_0) = O(\|x_0\| + F_0), \forall t \geq t_0 \text{ and } \lim_{t_0 \leq t \to \infty} \|Y_m(t, x_0)\| \leq O(F_0), \text{ which}$$

implies that under conditions of this statement $Y_m \in \Omega$.

Next, to embrace the application of (5.2) to (6.1), we will show by contradiction that, $z + Y_m \in \Omega$, $\forall t \geq t_0$ if $\|x_0\|$ and $F_0$ are sufficiently small. In fact, since $z(0) = 0$, let us assume that, $\exists t_1 > t_0$ such that, $z(t_1) + y(t_1) \in \partial(\Omega)$, where $\partial(\Omega)$ is a boundary of $\Omega$. Then, due to (6.6),

$$Z_1(t, x_0) \leq \int_{t_0}^{t_1} e^{-\hat{\lambda}(t-\tau)} \bar{\Gamma}(\tau, \|x_0\|, F_0) d\tau \leq \sup_{t_0 \leq t \leq t_1} \bar{\Gamma} \int_{t_0}^{t_1} e^{-\hat{\lambda}(t-\tau)} d\tau = O(\|x_0\| + F_0), \forall t \in [t_0, t_1].$$

This implies that $Z_1(t, x_0)$ can be made arbitrary small $\forall t \in [t_0, t_1]$ if $\|x_0\|$ and $F_0$ are sufficiently small which infers that, $z(t_1) + y(t_1) \in \Omega - \partial\Omega$ under the assumptions of this statement.

In addition, (6.6) implies that,

$$Z_1(t, x_0) \leq \bar{Z}_1(t, \|x_0\|, F_0) = \int_{t_0}^{t} e^{-\hat{\lambda}(t-\tau)} \bar{\Gamma}(\tau, \|x_0\|, F_0) d\tau = O(\|x_0\| + F_0) \int_{t_0}^{t} e^{-\hat{\lambda}(t-\tau)} d\tau \in \Theta_2, t \geq t_0; \text{ hence,}$$

$\bar{Z}_1(t, \|x_0\|, F_0)$ increases in both $\|x_0\|$ and $F_0$ for $\forall t \geq t_0$. In turn, $\|z(t, x_0)\| \leq \|z_1(t, x_0)\| \leq \|Z_1(t, x_0)\| < \infty$

, $t \geq t_0$ and $\lim_{t_0 \leq t \to \infty} \|z(t, x_0)\| = O(F_0)$. Consequently, $\|x(t, x_0)\| < \infty$, $\lim_{t_0 \leq t \to \infty} \|x(t, x_0)\| = O(F_0)$ and

$$\|Y_m(t, x_0)\| - \bar{Z}_1(t, \|x_0\|, F_0) \leq \|Y_m(t, x_0)\| - Z_1(t, x_0) \leq \|x(t, x_0)\| \leq \|Y_m(t, x_0)\| + Z_1(t, x_0)$$
$$\leq \bar{Y}_m(t, \|x_0\|, F_0) + \bar{Z}_1(t, \|x_0\|, F_0), \forall t \geq t_0 \ \square$$

**Theorem 8.** Assume that $F_0 = 0$, $\|x_0\|$ is sufficiently small, $f'(t, x)$ is a continuous matric, $\lambda(t) < -\hat{\lambda}, \forall t \geq t_0$, $\hat{\lambda} > 0$, inequalities (4.1), (5.2) and (5.3) hold and $c(t) \leq \hat{c}$. Then, the trivial solution to (2) is asymptotically stable and

$$\|Y_m(t, x_0)\| - \hat{Z}_1(t, \|x_0\|) \leq \|Y_m(t, x_0)\| - Z(t, x_0) \leq \|x(t, x_0)\| \leq \|Y_m(t, x_0)\| + Z_1(t, x_0)$$
$$\leq \hat{Y}_m(t, \|x_0\|) + \hat{Z}_1(t, \|x_0\|), \forall t \geq t_0$$
(6.8)

where $\hat{Y}_m(t, \|x_0\|) = \sum_{i=1}^{m} \hat{y}_i(t, \|x_0\|)$, $\hat{Z}_1(t, \|x_0\|) = \int_{t_0}^{t} e^{-\hat{\lambda}(t-\tau)} \hat{\Gamma}(\tau, \|x_0\|) d\tau$ and both left – side inequalities should be adjusted to zero if they assume negative values.

**Proof.** In fact, in this case, due to Theorem 1, $\|Y_m(t, x_0)\| \leq O(\|x_0\|), \forall t \geq t_0$ and $\lim_{t_0 \leq t \to \infty} \|Y_m(t, x_0)\| = 0$.

As in the prior statement, we can readily show that both $\sum_{i=1}^{m-1} y_i \in \Omega_1$ and $z + \sum_{i=1}^{m} y_i \in \Omega$, $\forall t \geq t_0$ for sufficiently small $\|x_0\|$. Hence, $Z_1(t, x_0) \leq \hat{Z}_1(t, \|x_0\|) = \int_{t_0}^{t} e^{-\hat{\lambda}(\tau - t_0)} \hat{\Gamma}(\tau, \|x_0\|) d\tau \in \Theta_1$ since, $\hat{\Gamma}(t, x_0) \in \Theta_1$. This implies that, $\lim_{t_0 \leq t \to \infty} \|z(t, x_0)\| = 0$ and that (6.8) holds $\square$

Note that in (6.7) and (6.8) most conservative upper bounds depend only upon scalar parameters, $\|x_0\|$ and $F_0$ or $\|x_0\|$, respectively, whereas in most conservative lower bounds just the estimate of error function depend upon these scalar parameters.

Let us assume, in turn, that the approximate solution $Y_m(t, x_0)$ is defined by (4.3). Clearly, in this case equation (6.3) holds with the same $\lambda(t, t_0)$ and formulas (6.2), (6.4) and (6.5) hold as well with obvious simplifications. Thus, theorems 7 and 8 hold in this case as well.

## 7. Nonlinear Error Equation and Its Applications

Utility of linear error equation (6.3) enables derivation of local stability and boundedness criteria but fails to estimate trapping/stability regions of the original equations. This section develops a nonlinear version of the error equation via application of a nonlinear extension of Lipschitz inequality, see [43] and section (2.1) of this paper for a brief introduction of this concept. Such nonlinear error equation frequently better reflects a nonlinear character of the underlying equations and enables estimation of the trapping/stability regions for systems (1) or (2), respectively.

### 7.1. Nonlinear Error Equation

To streamline further derivations, we assume in this section that $f(t, x)$ is a vector polynomial in $x$ without linear term and $f(t, 0) = 0$. Under these conditions, (4.2) can be written as follows,

$$\dot{z} = A(t)z + \pi(t, y, z)$$
$$z(t_0, x_0) = 0 \qquad (7.1)$$

where $y = [y_1, ..., y_m]^T$, $\pi(t, y, z)$ is a polynomial in $z$ with coefficients depending upon $t$ and $y$. Note that in the reminder of this section we assume that $y(t, x_0)$ is defined by either (4.1) or (4.3) under conditions of Theorem 4/Theorem 5/Theorem 6, respectively. Thus, $\|y(t, x_0)\| \leq \|y(t, x_0)\|_1 \leq \hat{Y}_m(t, \|x_0\|) = \sum_{i=1}^m \hat{y}_i(t, \|x_0\|) \in \Theta_1$ if $F_0 = 0$ and $\|y(t, x_0)\| \leq \overline{Y}_m(t, \|x_0\|, F_0) = \sum_{i=1}^m \overline{y}_i(t, \|x_0\|, F_0) \in \Theta_2$ if $F_0 > 0$. In turn, successive application of standard norm inequalities (see section 2.3 and Appendix of this paper) yields that,

$$\|\pi(t, y, z)\| \leq \Pi_*(t, y, \|z\|) + \gamma_*(t, y) \qquad (7.2)$$

where $\Pi_*(t, y, \|z\|) = \Pi(t, x_0, \|z\|)$ is a scalar-polynomial in $\|z\|$, $\Pi(t, x_0, 0) = 0$, $\|\pi(t, y, 0)\| = \gamma_*(t, y)$ $= \gamma(t, x_0) \in \mathbb{R}$ and $\gamma_*(t, y)$ is a scalar polynomial in $y$ with $\gamma_*(t, 0) = 0$. Additionally, $\gamma_*(t, y) \leq \hat{\gamma}_*(t, \hat{Y}_m)$ $= \hat{\gamma}(t, \|x_0\|) \in \Theta_1$ if $F_0 = 0$ and $\gamma_*(t, y) \leq \overline{\gamma}_*(t, \overline{Y}_m) = \overline{\gamma}(t, \|x_0\|, F_0) \in \Theta_2$ if $F_0 > 0$. Thus, $\lim_{t_0 \leq t \to \infty} \gamma(t, x_0) = 0$ if $F_0 = 0$ and $\lim_{t_0 \leq t \to \infty} \gamma(t, x_0) = O(F_0)$ if $F_0 > 0$, see Appendix for additional details on definitions of $\hat{\gamma}_*$ and $\overline{\gamma}_*$.

Consequently, application of (2.12) to (7.1) with account of (7.2) returns the following inequality,
$$D^+ \|z_2\| \leq p(t)\|z_2\| + c(t)\Pi(t, x_0, \|z_2\|) + c(t)\gamma(t, x_0), \forall t \geq t_0$$
$$\|z_2(t_0, x_0)\| = 0 \qquad (7.3)$$

Clearly, (7.3) implies that $\|z\| \leq \|z_2\|$.

Subsequent applications of the comparison principle to (7.3) yields the following scalar equation,
$$\dot{Z}_2 = p(t)Z_2 + c(t)\Pi(t, x_0, Z_2) + c(t)\gamma(t, x_0), \forall t \geq t_0$$
$$\|Z_2(t_0, x_0)\| = 0 \qquad (7.4)$$

where $\|z\| \leq \|z_2\| \leq Z_2$. In the reminder of this paper we assume that (7.4) admit a unique solution for $t \geq t_0$.

**7.2. Estimation of Bilateral Solutions' Norms and Trapping/Stability Regions**

The following statements infer some properties of solutions to multidimensional equations (1) and (2) under assumptions set on solutions to scalar equations (7.4). We assume in these statements without repetition that $\Omega_\times(t_0) \subset \mathbb{R}^n$ is a bounded neighborhood of $x \equiv 0$.

**Theorem 9**. Assume that $f(t, x)$ is a vector polynomial in $x$ without linear component, $f(t, 0) = 0$, $F_0 > 0$ is a sufficiently small value, $v_1 > 0$, $y_k(t, x_0)$, $k \in \mathbb{N}_m$ are solutions to (4.1), $\lim_{t_0 \leq t \to \infty} Z_2(t, x_0) = O(F_0)$, $\forall x_0 \in \Omega_\times(t_0)$ and $x(t, x_0)$ is a solution to (1).

Then, $\|x(t, x_0)\| < \infty$, $\forall t \geq t_0$, $\lim_{t_0 \leq t \to \infty} x(t, x_0) = O(F_0)$, $\forall x_0 \in \Omega_\times(t_0)$ and $\Omega_\times$ is included in the trapping region of (1). Furthermore, under the above conditions, solutions to (1) assume the following bilateral bounds,

$$\|Y_m(t, x_0)\| - Z_2(t, x_0) \leq \|x(t, x_0)\| \leq \|Y_m(t, x_0)\| + Z_2(t, x_0), \forall t \geq t_0 \qquad (7.5)$$

where the left side of (7.5) should be adjusted to zero if it takes negative values.

**Proof.** Really, in this case due to Theorem 5, $\|Y_m\| \leq \bar{Y}_m \in \Theta_2$, $t \geq t_0$ and, in turn, $\lim_{t_0 \leq t \to \infty} \|Y_m(t, x_0)\| \leq \lim_{t_0 \leq t \to \infty} \bar{Y}_m(t, \|x_0\|, F_0) = O(F_0)$. Hence, under conditions of this statement, (7.5) follows from (1.10) which, subsequently, implies that $x(t, x_0) < \infty$, $\forall x_0 \in \Omega_\times(t_0)$, $t \geq t_0$ and $\lim_{t_0 \leq t \to \infty} x(t, x_0) = O(F_0)$, $\forall x_0 \in \Omega_\times(t_0)$ □

**Theorem 10.** Assume that $f(t,x)$ is a vector polynomial in $x$ without linear component, $f(t,0) = 0$, $F = 0$, $v_1 > 0$, $y_k(t, x_0)$, $k \in \mathbb{N}_m$ are solutions to (4.1), $x(t, x_0)$ is a solution to (2) and $\lim_{t_0 \leq t \to \infty} Z_2(t, x_0) = 0$, $\forall x_0 \in \Omega_\times(t_0)$.

Then, the trivial solution to (1.2) is asymptotically stable, $\Omega_\times(t_0)$ is included in the stability region of the trivial solution to (1.2) and $\|x(t, x_0)\|$ admits bilateral bounds (7.5) with $F = 0$.

**Proof.** In fact, in this case, due to Theorem 4, $\lim_{t_0 \leq t \to \infty} \|Y_m(t, x_0)\| = 0$, $\forall x_0 \in \Omega_\times(t_0)$, hence, the above theorem follows from (1.10) and inequalities (7.5) with $F_0 = 0$ □

**Theorem 11.** Assume that $f(t, x)$ is a vector polynomial in $x$ without linear component, $f(t,0) = 0$, $v_1 > 0$, $y_k(t, x_0)$, $k \in \mathbb{N}_m$ are solutions to (4.3). Assume also that: 1. $\lim_{t_0 \leq t \to \infty} Z_2(t, x_0) = 0$, $\forall x_0 \in \Omega_\times(t_0)$ if $F_0 = 0$ and in this case $x(t, x_0)$ is a solution to (2); 2. $\lim_{t \to \infty} Z_2(t, x_0) = O(F_0)$, $\forall x_0 \in \Omega_\times(t_0)$ if $F_0 > 0$ and in this case $x(t, x_0)$ is a solution to (1).

Then, 1. the trivial solution to (1.2) is asymptotically stable, $\Omega_\times(t_0)$ is included in the stability region of the trivial solution to (1.2) and $\|x(t, x_0)\|$ admits bilateral bounds (7.5) with $F = 0$. 2. $\|x(t, x_0)\| < \infty$, $\forall t \geq t_0$, $\lim_{t_0 \leq t \to \infty} x(t, x_0) = O(F_0)$ and $\Omega_\times(t_0)$ is included in the trapping region of (1) containing $x \equiv 0$ and $\|x(t, x_0)\|$ admits bilateral bounds (7.5).

**Proof.** The proof of this statement is based on Theorem 6 and is analogous to the proofs of the previous two statements.

### 7.3 Approximate Solution to Error Equation

Though (7.4) is non-integrable scalar equation, the second and third additions in the right-hand side of this equation are nonnegative, whereas the first - can be either positive/negative or change sine for some values of $t$. To exploit this property, we set that $\Pi(t, x_0, Z_2) = D(t, x_0) Z_2 + \Pi_-(t, x_0, Z_2)$, where $0 \leq D(t, x_0) \in \mathbb{R}$, $t \geq t_0$, $\Pi_-(t, x_0, Z_2)$ is a scalar polynomial in $Z_2$ without linear component and with $\Pi_-(t, x_0, 0) = 0$. In turn, setting $\Pi_-(t, x_0, Z_2) = 0$ yields a linear, nonhomogeneous and scalar equation,

$$\dot{Z}_3 = P_1(t, \|x_0\|) Z_3 + c(t) \gamma(t, x_0), \quad \forall t \geq t_0$$
$$\|Z_3(t_0)\| = 0$$
(7.6)

where $P_1(t, x_0) = p(t) + c(t) D(t, x_0)$. Clearly, this equation possesses a unique solution

$$Z_3(t, x_0) = \int_{t_0}^{t} \theta_1(t, \tau, x_0) c(\tau) \gamma(\tau, x_0) d\tau \qquad (7.7)$$

where $\theta_1(t,\tau,x_0) = \exp\left(\int_\tau^t P_1(s,x_0)ds\right)$. Apparently, $Z_3(t,x_0) \leq Z_2(t,x_0)$, $t \geq t_0$, and $Z_2(t,x_0)$ $Z_3(t,x_0)$ should be close to each other if $\Pi_-(t,x_0,Z_2)$ is sufficiently small. Clearly, if $\lim_{t_0 \leq t \to \infty} Z_3(t,x_0) = \infty$, $\forall x_0 \in \omega(t_0)$, then $\lim_{t_0 \leq t \to \infty} Z_2(t,x_0) = \infty$, $\forall x_0 \in \omega(t_0)$ as well. Thus, $\Omega_\times(t_0) \subset \omega(t_0)$, where $\Omega_\times$ is an estimate of the trapping/stability region attained through evaluation of (7.4). Note also that $\Omega_\times(t_0)$ is enclosed in actual the trapping/stability region, whereas $\omega(t_0)$ - ought not be included in such region. Yet, our simulations show that $\omega(t_0)$ provides a valuable assessment of the actual trapping/stability regions which can be readily obtained by utility of (7.7).

Finally, we present a sufficient condition under which $\lim_{t_0 \leq t \to \infty} Z_3(t,x_0) = 0$

**Theorem 12.** Assume that $\lim_{t_0 < t \to \infty} \sup (t-t_0)^{-1} \int_{t_0}^t P_1(s,x_0)ds = -\lambda_1(t_0,x_0) < 0$, $\lambda_1 > 0$, $x_0 \in \omega(t_0) \subset \mathbb{R}^n$, $c(t) \leq \bar{c}$, $t \geq t_0$ and $y_k(t,x_0)$, $k \in \mathbb{N}_m$ are defined by either (4.1) or (4.3). Then, 1. $\lim_{t_0 \leq t \to \infty} Z_3(t,x_0) = 0$, $\forall x_0 \in \omega(t_0)$ if $F_0 = 0$ and 2. $\lim_{t_0 \leq t \to \infty} Z_3(t,x_0) = O(F_0)$, $\forall x_0 \in \omega(t_0)$ if $F_0 > 0$.

**Proof.** In fact, $\lambda_1(x_0)$ defines the Lyapunov exponent for the homogeneous part of (7.6). Hence, $\theta_1(t,\tau,x_0) = D_{1\varepsilon} \exp((\lambda_1+\varepsilon)(t-\tau))$, where constants $D_{1\varepsilon}$ and $\varepsilon > 0$, see section 2.1 of this paper and further citation therein. In turn, $Z_3(t,x_0) \leq D_{1\varepsilon} \bar{c} \int_{t_0}^t \exp((\lambda_1+\varepsilon)(t-\tau))\gamma(\tau,x_0)d\tau$. Thus,

$$\int_{t_0}^t \exp((\lambda_1+\varepsilon)(t-\tau))\gamma(\tau,x_0)d\tau \leq \int_{t_0}^t \exp((\lambda_1+\varepsilon)(t-\tau))\hat{\gamma}(\tau,\|x_0\|)d\tau \in \Theta_1, t \geq t_0 \text{ if } F_0 = 0 \text{ and}$$

$$\int_{t_0}^t \exp((\lambda_1+\varepsilon)(t-\tau))\gamma(\tau,x_0)d\tau \leq \int_{t_0}^t \exp((\lambda_1+\varepsilon)(t-\tau))\bar{\gamma}(\tau,\|x_0\|,F_0)d\tau \in \Theta_2, t \geq t_0 \text{ if } F_0 > 0 \; \square$$

## 8. Simulations

This section applies the developed above methodology to estimate the bounds of solutions and the trapping/stability regions for the Van der-Pol and Duffing -like models featuring typical dissipative and conservative nonlinearities.

### 8.1 Van der Pol – like Model

Firstly, we simulate system (3.1) with the following values of parameters which includes a large nonlinear coefficient, $\omega = 2$, $\alpha_1 = 1.2$, $\alpha_2 = -100$, $a_1 = a_2 = 5$, $r_1 = 4.8\pi$, $r_2 = 21$. Note that $a$ equals either zero or 0.23 in this set of simulations. For system (3.1) with $a = 0$ and a stable trivial solution typical time-histories of upper estimates of the norms of error terms are plotted in Figure 1 for initial vectors locating either within or beyond its stability region.

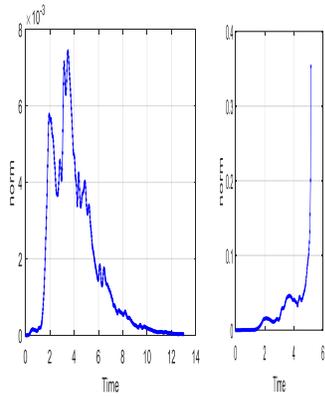
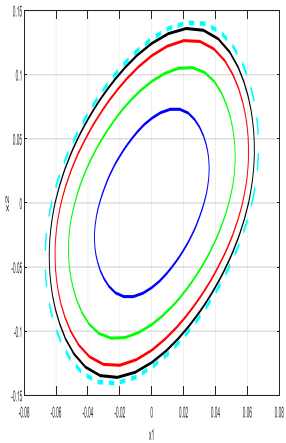
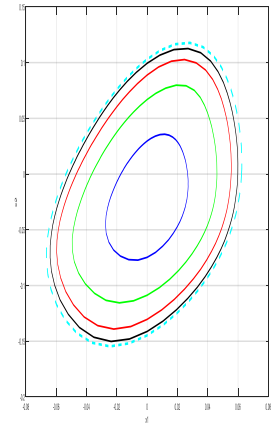

Fig. 2a Estimation of stability region of the trivial solution of (3.1) by the first approach

Fig. 2b Estimation of trapping region for (3.1) by the first approach

Fig. 1a    Fig. 1b

Fig. 1 Time histories of the estimates of the norm of error term computed for initial values within (Figure 1a) and beyond (Figure 1b) stability region of homogeneous counterpart of (3.1)

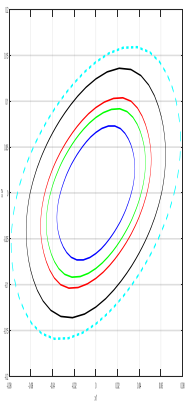
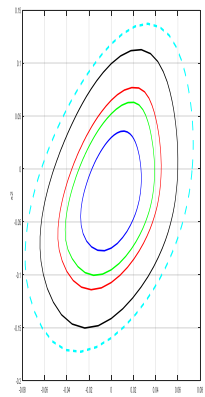
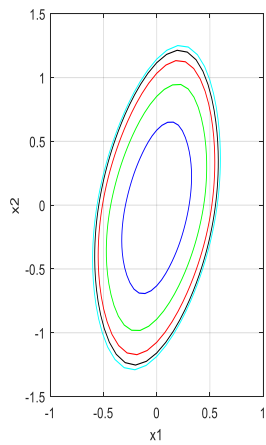
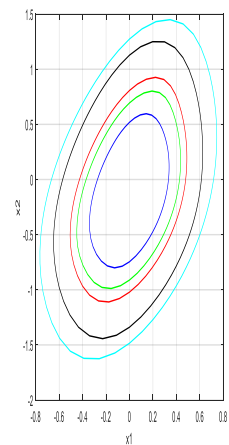

Fig. 2c Estimation of stability region of the trivial solution of (3.1) by second approach

Fig. 2d. Estimation of trapping region for (3.1) delivered by second approach.

Fig. 3a Estimation of a larger size stability regions for system (3.1) by first approach

Fig. 3b Estimation of a larger size stability regions for (3.1) by second approach

In both cases the estimates are stemmed from zero and approach either zero (Figure 1a) or infinity (Figure 1b) with $t \to \infty$. If $a \neq 0$ and $x_0$ is fitted into the trapping region of (3.1), the estimates approach $O(a)$; if in this case $x_0$ is located beyond the trapping region, then the estimate is approach infinity in our simulations. The threshold sets of $x_0$ differentiating these patterns of behavior of solutions to a scalar equation (7.4) yield the estimate of the boundaries of trapping/stability regions of (3.1) in our simulations.

Fig. 2 compare three estimates to the boundaries of the corresponding regions which are developed consecutively by simulating equations (3.2) and (3.4) (blue – lines), (3.5) and (3.7) (green -lines) and (3.8) and (3.10) (red - lines). As is mentioned in section 3, the complete set of equations in the last two cases include also the appropriate equations derived at the previous steps. Black – lines mark the reference estimates of actual boundaries of the corresponding regions, which are obtained in simulations of (3.1). Aqua – lines mark approximations delivered by linearized equation (3.10). These lines are located beyond the regions encircled by red – lines.

Figure 2a and 2b display reference and approximate boundaries of trapping/stability regions for (3.1)

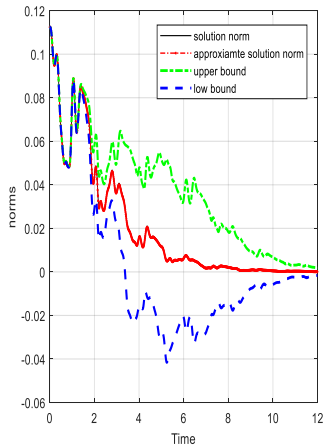

Fig. 4a displays simulated time – histories of the accurate and approximate solutions to homogeneous counterpart of (3.1) as well as the running low and upper bounds of the accurate solution stemmed from initial vectors located within stability region of (3.1) near its boundary

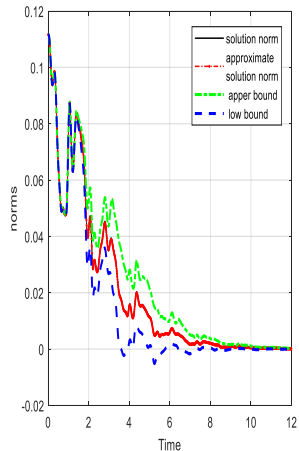

Fig. 4b displays the simulated time – histories of the accurate and approximate solutions and low and upper bounds of the accurate solution stemmed from the initial vector located in the central part of stability region of homogeneous counterpart of (3.1)

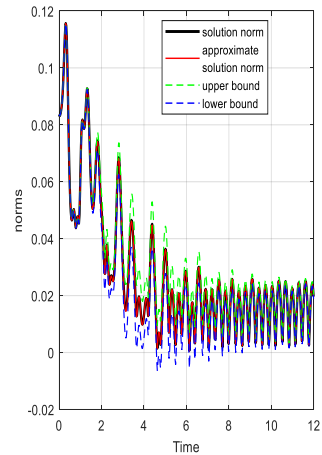

Fig. 4c displays simulated time – histories of the accurate and approximate solutions and running upper and lower bounds of the accurate solution stemmed from initial vector located within the trapping region of (3.1).

corresponding to our first approach. Figures 2c and 2d display in matching colors the corresponding approximations delivered by the second approach that is defined by equations (3.2), (3.11), (3.12) and (3.14). Note that in Figures 2a and 2c, $a = 0$ and in Figures 2b and 2d, $a = 0.23$. Clearly, successive approximations considerably enhance the accuracy of the corresponding estimates on every consecutive step for both approaches. More efficient in computations second approach provides less accurate estimates than its more elaborate counterpart for every used iteration. In the first approach, solutions to linearized and integrable error equations (aqua – lines) provides the estimates that are sufficiently close to the reference ones, but the differences between these estimates increase for less accurate second approach. Note that in all these cases the linearized equations provide excessive estimates of the trapping/stability regions.

Figure 3 shows application of the first approach to estimation of a larger size stability region of equation (3.1). In these simulations, $\alpha_1 = 1$, $\alpha_2 = -1$, $a = 0$ and other values of parameters are equal to the previous ones. Literally, the accuracy of our approximations remains intact in these cases as well.

Simulated time -histories of the norms of accurate and approximate solutions to (3.1) as well as their lower and upper bounds are displayed in Figure 4. Figures 4a and 4b show the results of simulations of (3.1) with $a = 0$ whereas Figure 4c display simulations made with $a = 0.23$. In all these three cases the difference between the norms of approximate and authentic solutions remain small if the initial vectors are located within the major parts of the trapping/stability regions but near its boundary. Consecutively, black and red – lines, which display simulated time – histories of the norms of accurate and approximate solutions, are practically overlaid on these figures. The running upper and lower bounds of the estimates of norms of accurate solutions are plotted in dashed green and blue lines.

Note that the lower bounds should be adjusted to zero if they take negative values. In fact, these estimates turn out to be less accurate if the initial vector is located near the region's boundary, see Figure 4a, but the accuracy provided by these bounds increases substantially if the initial vector is moved into a central part of stability region, see Figure 4b. Figure 4c plot the corresponding bounds and the norms of accurate and approximate solutions to (3.1) stemmed

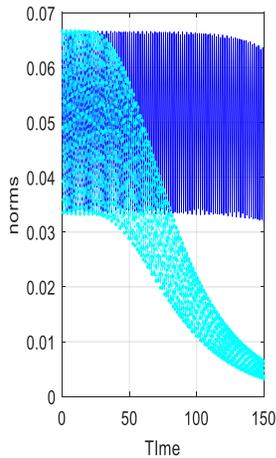 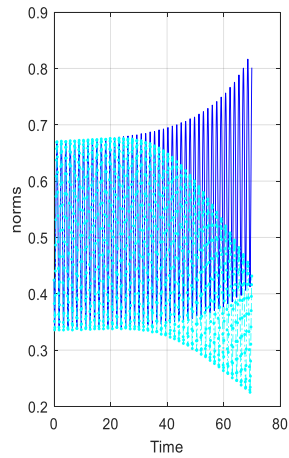 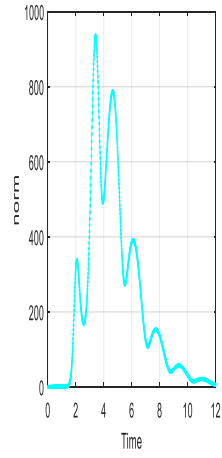 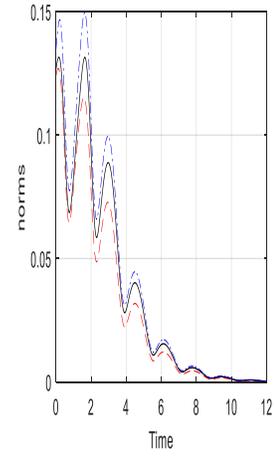

Fig. 5a displays in blue and aqua - lines simulated time – histories of the accurate and approximated solutions emanated from initial vectors locating within the stability region of homogeneous counterpart of (3.1) near its

Fig. 5b displays in blue and aqua - lines simulated time – histories of the accurate and approximated solutions emanated from initial vector locating beyond the stability region of homogeneous counterpart of (3.1) near its boundary

Fig. 5c displays simulated time – history of approximate solution emanated from initial vector locating beyond the stability region of homogeneous counterpart of (3.1)

Fig. 5d displays simulated time – histories of approximate solutions emanated from initial vectors locating within (red – line) and beyond (black and blue – lines) stability region of homogeneous counterpart of (3.1) near its from the initial vector locating within the trapping region of the "forced solution" to this equation. In all of these three cases time -histories of the estimates of the error – norms relatively quickly reach their maximal value and subsequently approach either zero or $O(a)$ if $t \to \infty$.

Analysis of behavior of the norms of approximate solutions in narrow neighborhoods of the boundaries of trapping/ stability regions spotlight some distinctiveness that can be used to estimate the boundaries of these regions via bypassing utility of (7.4). In fact, Figures 5a and 5b display behavior of the norms of accurate and approximated solutions emanated from initial vector locating within and beyond the stability region of homogeneous counterpart of (3.1) near its boundary. In both cases the norms of solutions are practically overlaid on initial time intervals but diverge thereafter. Figure 5c demonstrates that the maximal value of the norm of an approximate solution to this equation rapidly increases if the corresponding initial vector is moved beyond stability region. Finally, Figure 5d plots in red, black and blue – lines time -histories of the norms of approximate solutions stemmed from initial vectors located near stability boundary. The initial vector with smallest norm (red – line) is located within the stability region of this equation near its boundary, whereas two other initial vectors are located beyond the boundary of this region. In turn, first two local maxima are practically equal on the black – line, whereas on blue line, second maxima exceed the first one. Hence, an estimation of the boundary of the trapping/stability region can be rendered by selecting such initial vectors which practically equate first two local maximal values of time-histories of the corresponding norms of approximate solutions. Application of this heuristic criterion let us to skip simulation of (7.4) and in our simulations provides upper estimates of the norms of threshold initial vectors which exceed the norms of directly simulated ones by about 4% /5%.

### 8. 2. **Duffing – like Model**

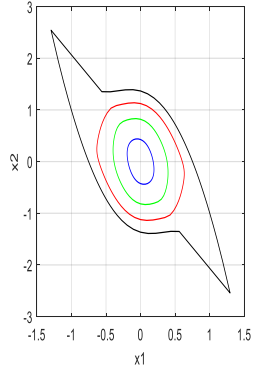

Fig. 6a plots the estimated boundaries of stability region for Duffing system

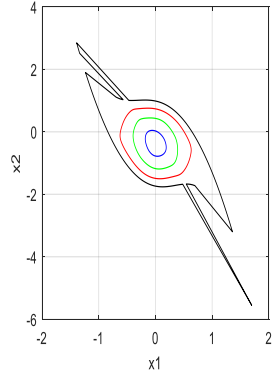

Fig. 6b plots the estimated boundaries of trapping region for Duffing system, $F_0 = 3.5$

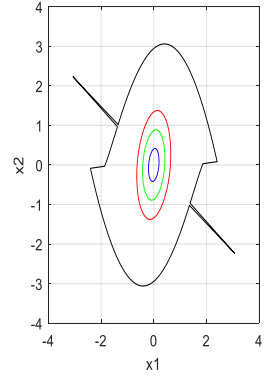

Fig. 6c plots the estimated boundaries of stability region for Duffing system with small damping

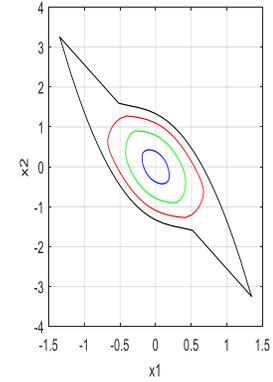

Fig. 6d plots the estimated boundaries of stability region for Duffing system, $t_0 = \pi / 2$

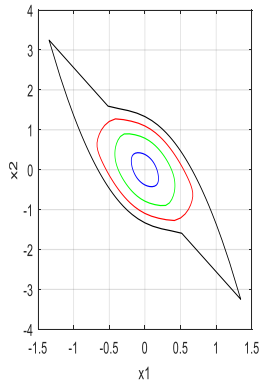

Fig. 6e plots the estimated boundaries of stability region for Duffing system, $t_0 = \pi^4$

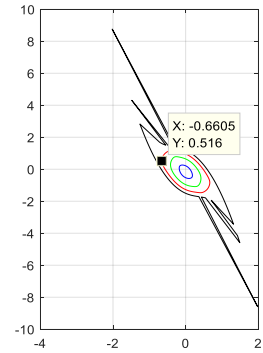

Fig. 6f plots the estimated boundaries of trapping region for Duffing system, $F_0 = 1.5$

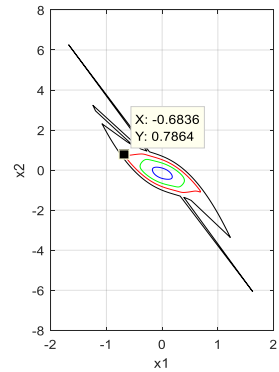

Fig. 6g plots the estimated boundaries of trapping region for Duffing system affected by initial pulse- function

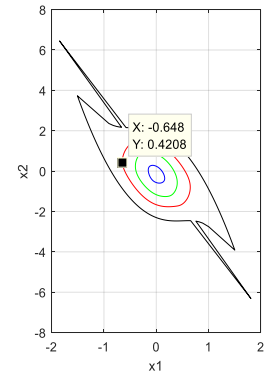

Fig. 6h plots the estimated boundaries of trapping region for Duffing system affected by initial pulse- function

This section applies the developed methodology for estimating the trapping/stability region of Duffing – like model which was widely used in various applications, see e.g. [47]. The governing equations for this model can be obtained by replacing the nonlinear term in the second equation of (3.1) as follows, $\alpha_2 x_2^3 \to \alpha_2 x_1^3$. Yet, the mechanisms involved in developing the trapping/stability regions in Van der-Pol and Duffing – like models are substantially different, but our methodology works in this case alike. Figure 6 display estimations of trapping/stability regions for Duffing - like models developed by our first approach which follows the protocol used prior in section 8.1. Figure 6a displays the estimates of stability boundaries for our model delivered by the first three approximations for the following values of parameters, $\alpha_1 = 1, \alpha_2 = -10, a = 0, a_1 = a_2 = 5$. Figure 6b displays the estimates of the trapping region of the forced solution of our model for $\alpha_1 = 1$, $\alpha_2 = -10$, $a = 3.5$, $a_1 = a_2 = 5$. Figure 6c shows that the accuracy of our approximations to stability boundaries decreases for small values of $\alpha_1 = 0.05$. Such behavior is expected since for this system the successive approximations diverge if $\alpha_1 = 0$. Nonetheless, utility of both extra additions in (1.10) and the transformations escalating the speed of convergence of these series, increase the accuracy of our approximations in these cases as well [48].

Figures 6d-6h show the impact of the change in $t_0$ on the shape of stability/trapping regions in our system. Figures 6a, 6d and 6e are computed for the same set of parameters but in the former case $t_0 = 0$ whereas in latter two cases $t_0 = \pi/2$ and $\pi^4$, respectively. It turns out that the differences in the shapes of the corresponding stability regions are obvious but rather <u>moderate</u> and consistent with almost – periodic character of variable coefficients entering this model.

Figures 6f-6h display the impact of the change in $t_0$ on the shape of trapping regions of the forced solution of Duffing system. Note that in these figures the most left points on the matching curves are labeled for the sake of comparison and $a_1 = a_2 = 0$. In all these cases $F_0 = 1.5$ and in Figure 6f, $\omega^2 = 4$, whereas in Figures 6g and 6h, $\omega^2(t) = 4 + \omega_1(t)$, where $\omega_1(t) = \nu(1 - H(t - \pi/2))$, $H(t)$ is a Heaviside function and $\nu$ equals either -2 or 4 in Figures 6g and 6h, respectively. Apparently, the pulse does not affect behavior of this system's solutions if $t_0 > \pi/2$, see Figure 6f.

Clearly, in these cases the effect of alteration of $t_0$ on the shapes of trapping regions of the corresponding system is noticeable but rather moderate as well. Meanwhile, piece-wise constant and other types of time – varying coefficients can seemingly change the behavior of solutions on some finite/infinite time – intervals which, in turn, will affect the shapes of trapping/stability regions.

Note also that our estimates tend to smooth the actual boundaries of trapping/stability regions which might have highly irregular structures [21].

## 9. Conclusion and Future Work

This paper develops a novel methodology for successive approximation of solutions which are stemmed from the trapping/stability regions of some nonlinear and time – varying systems. Next, we estimate the norm of error - component developed by such approximations. This let us to infer enhanced stability/ boundedness criteria, develop bilateral bounds of the norms of solutions and estimate the trapping/stability regions for such systems.

Two techniques were devised for successive approximations of solutions to original systems: first - conveys more accurate approximations and second - is more efficient in computations. Subsequently, we provide two approaches to estimation of the norms of error - components which engage the methodology that we developed in [43]. The first approach exploits application of the Lipschitz continuity condition which leads to scalar, linear and integrable error equation. The second approach adopts a nonlinear extension of Lipschitz condition and develops a nonlinear scalar error -equation with variable coefficients that better reflects the character of the underlined systems and allows estimation of their trapping/stability regions. Both approaches provide bilateral bounds for the norms of solutions and boundedness/stability criteria.

We show in inclusive simulations that just a few iterations provided by developed techniques return the approximations which resemble with high fidelity the time- histories of the norms of accurate solutions stemming from the trapping/stability regions of original systems. The accuracy of these approximations remains intact for nonlinearities of different kinds and magnitudes but can be affected if the largest Lyapunov exponent of the linear block of the corresponding system turn out to be close to zero. Consequently, such approximations enable recursive estimation of the trapping/stability regions.

Subsequently, we intend to develop upper and lower bounds to solutions of the nonlinear error – equation in close form, extend our methodology to systems under uncertainty and mitigate some of the technical constrains of the developed approach. This should allow us to extend this methodology to a wider class of systems arising in various application domains.

## Appendix

Let us sketch, for instance, how to derive inequalities (7.2). Assume that, $f(t,x) = [f_1,...,f_n]^T$, $x, u, z \in \mathbb{R}^n$, $f_i(t,x) = \zeta_i(t) x_i^{d_i} = \zeta_i(t)(z_i + u_i)^{d_i}$, $d_i \in \mathbb{N}$. Thus, $|f_i| = \left|\zeta_i(z_i + u_i)^{d_i}\right| \leq |\zeta_i| |z_i + u_i|^{d_i}$
$\leq |\zeta_i|(|z_i| + |u_i|)^{d_i} \leq |\zeta_i|(\|z\| + |u_i|)^{d_i} \leq |\zeta_i|(\|z\| + \|u\|)^{d_i}$. Then, taking into account that, $\|f\| \leq \|f\|_1 = \sum_{i=1}^{n} |f_i|$, we infer that, $\gamma_*(t,u) = \sum_{i=1}^{n} |\zeta_i| |u_i|^{d_i}$, $\Pi_*(t,y) = \sum_{i=1}^{n} |\zeta_i|(\|z\| + |u_i|)^{d_i} - \gamma_*(t,u)$, $\hat{\gamma}(t, \|u\|) = \sum_{i=1}^{n} |\zeta_i| \|u\|^{d_i}$ if $F_0 = 0$. Subsequently, we enter in this formula that $\|u\| \leq \|u\|_1$. The prior formula holds if $F_0 > 0$ with obvious change in notation.

Assume next that, $f_{id}(t,x) = \zeta_{id}(t) x^d$, where $x^d = x_1^{d_1}...x_n^{d_n}$, $d_i \in \mathbb{Z}$, $d = \sum_{i=1}^{n} d_i$. Then, $|f_{id}(t,x)| = |\zeta_{id}| \left|x_1^{d_1}\right|...\left|x_n^{d_n}\right|$ which enables utility of the prior procedure for every $\left|x_i^{d_i}\right|$. Clearly, such sequence of steps yields (7.2) if $f(t,x)$ is a vector - polynomial in $x$.

**Acknowledgement.** The programs used to simulate the models in section 8 of this paper were developed by Steve Koblik who also was engaged in discussions of the corresponding results.

Conflict of Interest: The author declares that he has no conflict of interest